\documentclass[10pt]{amsart}
\usepackage{amsthm,amsmath,amsfonts,amssymb,graphicx,epsf,epstopdf,pifont,subfigure,flafter,bm,amscd}
\usepackage{overpic}
\usepackage[margin=3.6cm]{geometry}
\usepackage{multirow}
\usepackage{capt-of}
\usepackage{xcolor}
\usepackage{tikz}
\AddToHook{env/table/begin}{\setlength{\belowcaptionskip}{0pt}}

\newcommand{\R}{\mathbb{R}}
\newcommand{\vf}{\mathcal{X}}

\newtheorem{lm}{Lemma}[section]

\newtheorem{prop}{Proposition}[section]

\newtheorem{re}{Remark}[section]
\newtheorem{cor}{Corollary}[section]
\newtheorem{maintheorem}{Theorem}

\numberwithin{equation}{section}

\begin{document}

\title[]
{Quartic polynomial differential systems with three disjoint invariant circles}

\author[C. Cao]
{Chen Cao${}^{*}$}
\thanks{${}^{*}$Corresponding author.}

\address{\textsc{Chen Cao}\newline
School of Mathematical Sciences, CMA-Shanghai, Shanghai Jiao Tong University,
Shanghai, 200240, P. R. China}
\email{aceaccce@sjtu.edu.cn}

\author[J. Llibre]
{Jaume Llibre}

\address{\textsc{Jaume Llibre}\newline
Departament de Matem\`{a}tiques, Universitat Aut\`{o}noma de Barcelona,
08193 Bellaterra, Barcelona, Catalonia, Spain \newline
Reial Acad\`emia de Ci\`encies i Arts de Barcelona, La Rambla 115,
08002 Barcelona, Catalonia, Spain}
\email{jaume.llibre@uab.cat}



\begin{abstract}
We determine the normal forms of real planar polynomial differential
systems of degree at most four with three prescribed pairwise disjoint
invariant circles, covering all configurations of their centers.
None of the prescribed circles is a limit cycle in the concentric case
or in the coaxial case with three distinct centers, where the pairwise radical axes coincide.
When exactly two centers coincide, or when the three centers
are distinct and collinear and the circles are non-coaxial,
we give necessary and sufficient conditions for
all three circles to be algebraic limit cycles, which are necessarily hyperbolic.
For non-collinear centers, we give explicit necessary and sufficient
conditions for all three circles to be periodic orbits and for all three
to be hyperbolic algebraic limit cycles. Each of the four configurations
with non-collinear centers is realized by a coprime quartic polynomial differential system with
three hyperbolic circular limit cycles.
\end{abstract}

\maketitle

\section{Introduction and statement of the main results}

Algebraic limit cycles are isolated periodic orbits that are ovals of real
invariant algebraic curves. For an overview of their study in planar
polynomial differential systems see \cite{Llibre2025,LlibreZhang2021}.

Christopher \cite{Christopher2001} constructed, for every
nonsingular real algebraic curve of degree $m$, a polynomial vector field
of degree $m$ whose only limit cycles are the bounded components of the
curve, all of them hyperbolic. Applied to the product of three pairwise
disjoint circle polynomials, this construction gives a system of degree six.

Llibre and Rodr\'{i}guez \cite{LlibreRodriguez2004} proved that every finite
configuration of pairwise disjoint simple closed curves in the plane can
be realized, up to homeomorphism, as the set of limit cycles of a polynomial
vector field, with all these limit cycles algebraic. Margalef-Bentabol and
Peralta-Salas \cite{MargalefBentabolPeraltaSalas2016} extended this construction
to prescribe periods, multiplicities and compatible stabilities. These
realization results do not guarantee a vector field of a prescribed degree.

For circular algebraic limit cycles the maximum number in a polynomial
vector field of degree $n$ is $n-1$, and the bound is sharp
\cite[Theorem~3]{LlibreRamirez2016}.
The same work gives an explicit quartic polynomial
differential system with three circular algebraic limit cycles and
non-collinear centers.

The inverse problem of determining polynomial vector fields with prescribed
invariant algebraic curves was studied by Christopher, Llibre, Pantazi and
Walcher \cite{ChristopherLlibrePantaziWalcher2009}, who obtained explicit
representations under nondegeneracy conditions and developed computational
methods for more general curves.
Gin\'{e}, Llibre and Valls \cite{GineLlibreValls2018} studied cubic
systems with two invariant circles and characterized when both circles are
algebraic limit cycles. Llibre, Messias and Reinol
\cite{LlibreMessiasReinol2017} obtained normal forms for quadratic and cubic
systems with two nonconcentric invariant circles and a prescribed Darboux
first integral, and classified the global phase portraits of a subclass of these systems.

In this paper we classify all real polynomial differential
systems of degree at most four preserving three prescribed pairwise of disjoint
circles, and for every configuration of their centers, we characterize when
all three prescribed circles are hyperbolic limit cycles. For non-collinear
centers, we also give a necessary and sufficient condition for all
three circles to be periodic orbits.
No first integral is assumed, and the classification does not require
these to be the only invariant circles of the system.

The main difficulty is to enforce simultaneous invariance
of the three prescribed circles without exceeding degree four, while
accounting for the different algebraic relations arising from their center
configurations. Our approach combines explicit conditions for the
cancellation of higher-degree terms with direct evaluation of the integrals
governing hyperbolicity, keeping the criteria for periodicity and
hyperbolicity separate.

More precisely we consider real polynomial differential systems of the form
\begin{equation}\label{sys:general}
\dot x=P(x,y):=\sum_{0\le i+j\le4}a_{ij}x^iy^j,
 \qquad \dot y=Q(x,y):=\sum_{0\le i+j\le4}b_{ij}x^iy^j,
\end{equation}
where $a_{ij},b_{ij}\in\R$.
The polynomial vector field associated with system~\eqref{sys:general} is
\[
 \vf=P\frac{\partial}{\partial x}+Q\frac{\partial}{\partial y}.
\]
The degree of polynomial differential system~\eqref{sys:general} is $\max\{\deg P,\deg Q\}\le4$;
the system is quartic when this maximum equals four.

Let $f\in\R[x,y]$ be non-constant.  The real algebraic curve
$f=0$ is invariant under system~\eqref{sys:general} if there is a polynomial
$K\in\R[x,y]$, called the cofactor of $f=0$, such that
\begin{equation}\label{eq:invariance}
 \vf(f):=P\frac{\partial f}{\partial x}
      +Q\frac{\partial f}{\partial y}=Kf.
\end{equation}
Here $\deg K\le3$.
We write $f\doteq g$ for equality up to a non-zero real constant and
understand greatest common divisors in this sense; $\gcd(P,Q)=1$
means that $P$ and $Q$ have no non-constant common factor.
The system is called coprime if $\gcd(P,Q)=1$.

A non-constant $C^1$ function $H$ on an open set
$U\subset\mathbb R^2$ is a first integral if $\vf(H)=0$ on $U$,
or equivalently, if $H$ is constant along solution curves in $U$.

A periodic orbit $\Gamma$, parametrized by a solution $\gamma(t)$
of minimal period $T>0$, is called hyperbolic if
\[
 \int_0^T(P_x+Q_y)(\gamma(t))\,dt\ne0.
\]
Equivalent criteria and stability properties are given in
Section~\ref{subsec:characteristic-exponents}.

A period annulus is an open annular region foliated by periodic orbits.

Let the three prescribed circles be
\[
 f_i=(x-a_i)^2+(y-b_i)^2-r_i^2=0,\qquad
 C_i=(a_i,b_i),\qquad r_i>0,\qquad i=1,2,3,
\]
and write $D_i$ for the closed disk bounded by $f_i=0$.  Throughout the paper
the three circles are assumed to be pairwise disjoint: for every $i<j$ either
\[
 d_{ij}>r_i+r_j\qquad\hbox{or}\qquad
 d_{ij}<|r_i-r_j|,\qquad d_{ij}=|C_i-C_j|.
\]
These alternatives correspond to exterior and strictly nested
circles, respectively, and exclude intersections and tangencies.

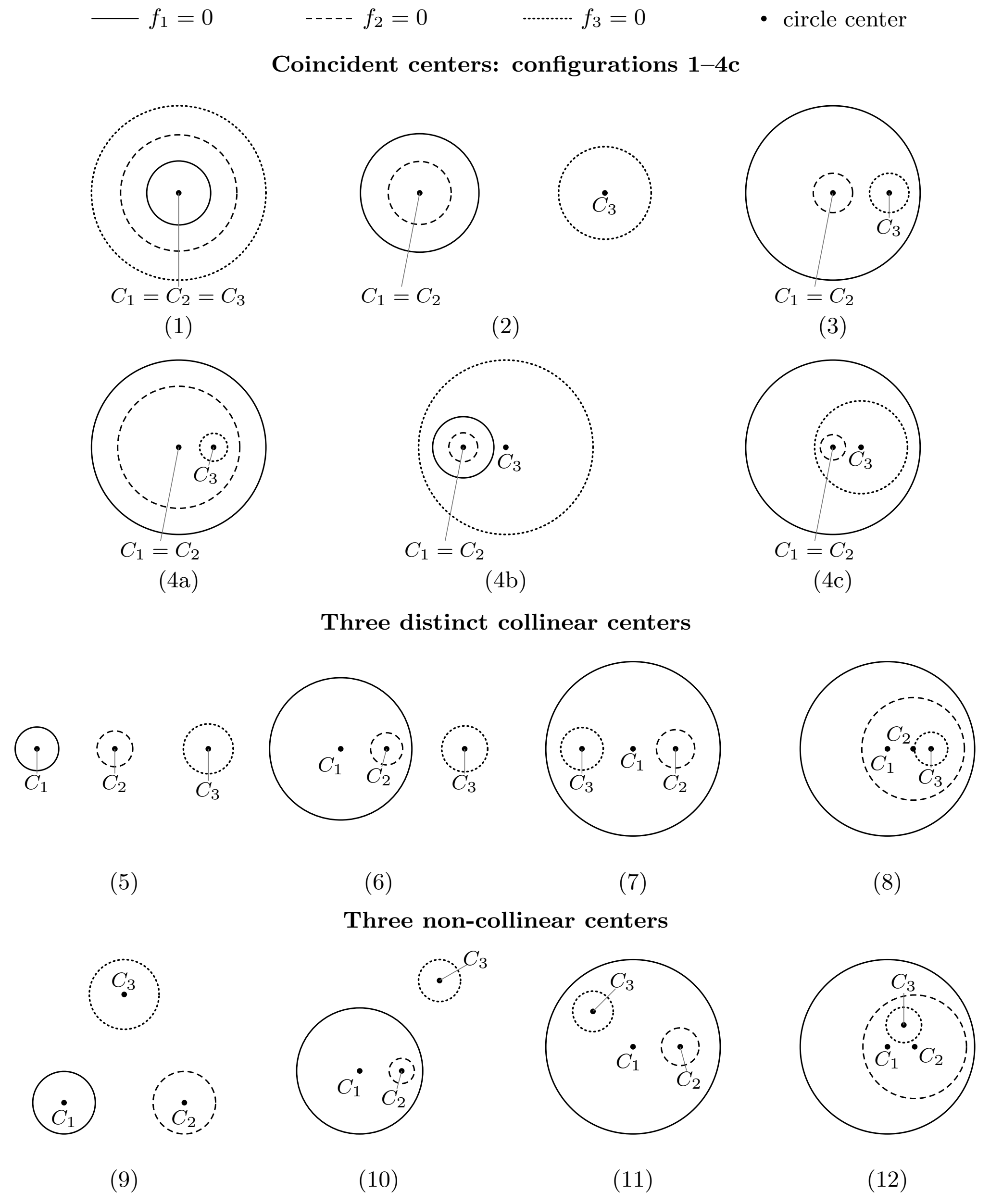
\begin{figure}[p]
\centering
\begin{tikzpicture}[x=1cm,y=1cm,
circleone/.style={black,line width=.55pt},
circletwo/.style={black,line width=.55pt,dash pattern=on 2.4pt off 1.5pt},
circlethree/.style={black,line width=.65pt,dash pattern=on .5pt off 1.3pt,line cap=round},
centerlabel/.style={font=\fontsize{8}{9.5}\selectfont,inner sep=1pt},
panelnumber/.style={font=\small}]
\path[use as bounding box] (0,-16.1) rectangle (14,0.55);
\draw[circleone] (1.3,0.25)--(1.95,0.25) node[right,font=\small] {$f_1=0$};
\draw[circletwo] (4.25,0.25)--(4.9,0.25) node[right,font=\small] {$f_2=0$};
\draw[circlethree] (7.25,0.25)--(7.9,0.25) node[right,font=\small] {$f_3=0$};
\fill (10.55,0.25) circle[radius=1.1pt] node[right=4pt,font=\small] {circle center};
\node[font=\small\bfseries] at (7,-.40) {Coincident centers: configurations 1--4c};
\node[font=\small\bfseries] at (7,-8.05) {Three distinct collinear centers};
\node[font=\small\bfseries] at (7,-12.15) {Three non-collinear centers};
\begin{scope}[shift={(2.50000,-2.15000)}]
\draw[circleone] (0.00000,0.00000) circle[radius=0.44000cm];
\draw[circletwo] (0.00000,0.00000) circle[radius=0.80000cm];
\draw[circlethree] (0.00000,0.00000) circle[radius=1.20000cm];
\fill (0.00000,0.00000) circle[radius=1.1pt];
\draw[black!55,line width=.3pt] (0.00000,0.00000) -- (0.00000,-1.29000);
\node[centerlabel] at (0.00000,-1.43000) {$C_1=C_2=C_3$};
\node[panelnumber] at (0,-1.85) {(1)};
\end{scope}
\begin{scope}[shift={(7.00000,-2.15000)}]
\draw[circleone] (-1.18471,0.00000) circle[radius=0.81529cm];
\draw[circletwo] (-1.18471,0.00000) circle[radius=0.43312cm];
\draw[circlethree] (1.36306,0.00000) circle[radius=0.63694cm];
\fill (-1.18471,0.00000) circle[radius=1.1pt];
\draw[black!55,line width=.3pt] (-1.18471,0.00000) -- (-1.43471,-1.29000);
\node[centerlabel] at (-1.43471,-1.43000) {$C_1=C_2$};
\fill (1.36306,0.00000) circle[radius=1.1pt];
\node[centerlabel] at (1.36306,-0.18000) {$C_3$};
\node[panelnumber] at (0,-1.85) {(2)};
\end{scope}
\begin{scope}[shift={(11.50000,-2.15000)}]
\draw[circleone] (0.00000,0.00000) circle[radius=1.20000cm];
\draw[circletwo] (0.00000,0.00000) circle[radius=0.27097cm];
\draw[circlethree] (0.77419,0.00000) circle[radius=0.27097cm];
\fill (0.00000,0.00000) circle[radius=1.1pt];
\draw[black!55,line width=.3pt] (0.00000,0.00000) -- (-0.25000,-1.29000);
\node[centerlabel] at (-0.25000,-1.43000) {$C_1=C_2$};
\fill (0.77419,0.00000) circle[radius=1.1pt];
\draw[black!55,line width=.3pt] (0.77419,0.00000) -- (0.77419,-0.34000);
\node[centerlabel] at (0.77419,-0.48000) {$C_3$};
\node[panelnumber] at (0,-1.85) {(3)};
\end{scope}
\begin{scope}[shift={(2.50000,-5.65000)}]
\draw[circleone] (0.00000,0.00000) circle[radius=1.20000cm];
\draw[circletwo] (0.00000,0.00000) circle[radius=0.84000cm];
\draw[circlethree] (0.48000,0.00000) circle[radius=0.19200cm];
\fill (0.00000,0.00000) circle[radius=1.1pt];
\draw[black!55,line width=.3pt] (0.00000,0.00000) -- (-0.25000,-1.29000);
\node[centerlabel] at (-0.25000,-1.43000) {$C_1=C_2$};
\fill (0.48000,0.00000) circle[radius=1.1pt];
\draw[black!55,line width=.3pt] (0.48000,0.00000) -- (0.41271,-0.25114);
\node[centerlabel] at (0.37647,-0.38637) {$C_3$};
\node[panelnumber] at (0,-1.85) {(4a)};
\end{scope}
\begin{scope}[shift={(7.00000,-5.65000)}]
\draw[circleone] (-0.58537,0.00000) circle[radius=0.42146cm];
\draw[circletwo] (-0.58537,0.00000) circle[radius=0.19902cm];
\draw[circlethree] (0.00000,0.00000) circle[radius=1.20000cm];
\fill (-0.58537,0.00000) circle[radius=1.1pt];
\draw[black!55,line width=.3pt] (-0.58537,0.00000) -- (-0.83537,-1.29000);
\node[centerlabel] at (-0.83537,-1.43000) {$C_1=C_2$};
\fill (0.00000,0.00000) circle[radius=1.1pt];
\node[centerlabel] at (0.05694,-0.21250) {$C_3$};
\node[panelnumber] at (0,-1.85) {(4b)};
\end{scope}
\begin{scope}[shift={(11.50000,-5.65000)}]
\draw[circleone] (0.00000,0.00000) circle[radius=1.20000cm];
\draw[circletwo] (0.00000,0.00000) circle[radius=0.17419cm];
\draw[circlethree] (0.38710,0.00000) circle[radius=0.63871cm];
\fill (0.00000,0.00000) circle[radius=1.1pt];
\draw[black!55,line width=.3pt] (0.00000,0.00000) -- (-0.25000,-1.29000);
\node[centerlabel] at (-0.25000,-1.43000) {$C_1=C_2$};
\fill (0.38710,0.00000) circle[radius=1.1pt];
\node[centerlabel] at (0.38710,-0.18000) {$C_3$};
\node[panelnumber] at (0,-1.85) {(4c)};
\end{scope}
\begin{scope}[shift={(1.75000,-9.80000)}]
\draw[circleone] (-1.20000,0.00000) circle[radius=0.30000cm];
\draw[circletwo] (-0.12857,0.00000) circle[radius=0.24643cm];
\draw[circlethree] (1.15714,0.00000) circle[radius=0.34286cm];
\fill (-1.20000,0.00000) circle[radius=1.1pt];
\draw[black!55,line width=.3pt] (-1.20000,0.00000) -- (-1.20000,-0.34000);
\node[centerlabel] at (-1.20000,-0.48000) {$C_1$};
\fill (-0.12857,0.00000) circle[radius=1.1pt];
\draw[black!55,line width=.3pt] (-0.12857,0.00000) -- (-0.12857,-0.34000);
\node[centerlabel] at (-0.12857,-0.48000) {$C_2$};
\fill (1.15714,0.00000) circle[radius=1.1pt];
\draw[black!55,line width=.3pt] (1.15714,0.00000) -- (1.15714,-0.44000);
\node[centerlabel] at (1.15714,-0.58000) {$C_3$};
\node[panelnumber] at (0,-1.85) {(5)};
\end{scope}
\begin{scope}[shift={(5.25000,-9.80000)}]
\draw[circleone] (-0.52105,0.00000) circle[radius=0.97895cm];
\draw[circletwo] (0.11053,0.00000) circle[radius=0.22105cm];
\draw[circlethree] (1.18421,0.00000) circle[radius=0.31579cm];
\fill (-0.52105,0.00000) circle[radius=1.1pt];
\node[centerlabel] at (-0.65605,-0.23383) {$C_1$};
\fill (0.11053,0.00000) circle[radius=1.1pt];
\draw[black!55,line width=.3pt] (0.11053,0.00000) -- (0.04323,-0.25114);
\node[centerlabel] at (0.00700,-0.38637) {$C_2$};
\fill (1.18421,0.00000) circle[radius=1.1pt];
\draw[black!55,line width=.3pt] (1.18421,0.00000) -- (1.18421,-0.34000);
\node[centerlabel] at (1.18421,-0.48000) {$C_3$};
\node[panelnumber] at (0,-1.85) {(6)};
\end{scope}
\begin{scope}[shift={(8.75000,-9.80000)}]
\draw[circleone] (0.00000,0.00000) circle[radius=1.20000cm];
\draw[circletwo] (0.58537,0.00000) circle[radius=0.26341cm];
\draw[circlethree] (-0.70244,0.00000) circle[radius=0.29268cm];
\fill (0.00000,0.00000) circle[radius=1.1pt];
\node[centerlabel] at (-0.00000,-0.18000) {$C_1$};
\fill (0.58537,0.00000) circle[radius=1.1pt];
\draw[black!55,line width=.3pt] (0.58537,0.00000) -- (0.58537,-0.34000);
\node[centerlabel] at (0.58537,-0.48000) {$C_2$};
\fill (-0.70244,0.00000) circle[radius=1.1pt];
\draw[black!55,line width=.3pt] (-0.70244,0.00000) -- (-0.70244,-0.34000);
\node[centerlabel] at (-0.70244,-0.48000) {$C_3$};
\node[panelnumber] at (0,-1.85) {(7)};
\end{scope}
\begin{scope}[shift={(12.25000,-9.80000)}]
\draw[circleone] (0.00000,0.00000) circle[radius=1.20000cm];
\draw[circletwo] (0.35294,0.00000) circle[radius=0.70588cm];
\draw[circlethree] (0.60000,0.00000) circle[radius=0.22941cm];
\fill (0.00000,0.00000) circle[radius=1.1pt];
\node[centerlabel] at (-0.05694,-0.21250) {$C_1$};
\fill (0.35294,0.00000) circle[radius=1.1pt];
\node[centerlabel] at (0.16202,0.19092) {$C_2$};
\fill (0.60000,0.00000) circle[radius=1.1pt];
\draw[black!55,line width=.3pt] (0.60000,0.00000) -- (0.60000,-0.26000);
\node[centerlabel] at (0.60000,-0.40000) {$C_3$};
\node[panelnumber] at (0,-1.85) {(8)};
\end{scope}
\begin{scope}[shift={(1.75000,-13.90000)}]
\draw[circleone] (-0.82759,-0.76966) circle[radius=0.43034cm];
\draw[circletwo] (0.82759,-0.76966) circle[radius=0.43034cm];
\draw[circlethree] (0.00000,0.72000) circle[radius=0.48000cm];
\fill (-0.82759,-0.76966) circle[radius=1.1pt];
\node[centerlabel] at (-0.82759,-0.98966) {$C_1$};
\fill (0.82759,-0.76966) circle[radius=1.1pt];
\node[centerlabel] at (0.82759,-0.98966) {$C_2$};
\fill (0.00000,0.72000) circle[radius=1.1pt];
\node[centerlabel] at (0.00000,0.90000) {$C_3$};
\node[panelnumber] at (0,-1.85) {(9)};
\end{scope}
\begin{scope}[shift={(5.25000,-13.90000)}]
\draw[circleone] (-0.26024,-0.33253) circle[radius=0.86747cm];
\draw[circletwo] (0.31807,-0.33253) circle[radius=0.17349cm];
\draw[circlethree] (0.83855,0.91084) circle[radius=0.28916cm];
\fill (-0.26024,-0.33253) circle[radius=1.1pt];
\node[centerlabel] at (-0.39524,-0.56636) {$C_1$};
\fill (0.31807,-0.33253) circle[radius=1.1pt];
\draw[black!55,line width=.3pt] (0.31807,-0.33253) -- (0.25078,-0.58367);
\node[centerlabel] at (0.21454,-0.71890) {$C_2$};
\fill (0.83855,0.91084) circle[radius=1.1pt];
\draw[black!55,line width=.3pt] (0.83855,0.91084) -- (1.21961,1.13084);
\node[centerlabel] at (1.34085,1.20084) {$C_3$};
\node[panelnumber] at (0,-1.85) {(10)};
\end{scope}
\begin{scope}[shift={(8.75000,-13.90000)}]
\draw[circleone] (0.00000,0.00000) circle[radius=1.20000cm];
\draw[circletwo] (0.64865,0.00000) circle[radius=0.25946cm];
\draw[circlethree] (-0.55135,0.48649) circle[radius=0.27892cm];
\fill (0.00000,0.00000) circle[radius=1.1pt];
\node[centerlabel] at (-0.05694,-0.21250) {$C_1$};
\fill (0.64865,0.00000) circle[radius=1.1pt];
\draw[black!55,line width=.3pt] (0.64865,0.00000) -- (0.73665,-0.32841);
\node[centerlabel] at (0.77288,-0.46364) {$C_2$};
\fill (-0.55135,0.48649) circle[radius=1.1pt];
\draw[black!55,line width=.3pt] (-0.55135,0.48649) -- (-0.24022,0.79761);
\node[centerlabel] at (-0.14123,0.89661) {$C_3$};
\node[panelnumber] at (0,-1.85) {(11)};
\end{scope}
\begin{scope}[shift={(12.25000,-13.90000)}]
\draw[circleone] (0.00000,0.00000) circle[radius=1.20000cm];
\draw[circletwo] (0.37500,0.00000) circle[radius=0.71250cm];
\draw[circlethree] (0.22500,0.30000) circle[radius=0.24375cm];
\fill (0.00000,0.00000) circle[radius=1.1pt];
\node[centerlabel] at (-0.00000,-0.18000) {$C_1$};
\fill (0.37500,0.00000) circle[radius=1.1pt];
\node[centerlabel] at (0.60883,-0.13500) {$C_2$};
\fill (0.22500,0.30000) circle[radius=1.1pt];
\draw[black!55,line width=.3pt] (0.22500,0.30000) -- (0.22500,0.74000);
\node[centerlabel] at (0.22500,0.88000) {$C_3$};
\node[panelnumber] at (0,-1.85) {(12)};
\end{scope}
\end{tikzpicture}
\caption{Geometric configurations of three pairwise disjoint
circles. Dots mark the centers; coincident centers share a dot. Solid,
dashed and dotted curves represent $f_1=0$, $f_2=0$ and $f_3=0$,
respectively. In configurations 2--4c $f_1=0$ and $f_2=0$ are always
the concentric pair.}
\label{fig:three-circle-configurations}
\end{figure}

The geometry of the centers and the containment relations
among the circles yield the following
configurations, shown in Figure~\ref{fig:three-circle-configurations}.

\begin{enumerate}
\item[(I)] \emph{Three coincident centers.}
  \begin{enumerate}
  \item[1.] three concentric circles.
  \end{enumerate}
\item[(II)] \emph{Exactly two coincident centers.}
  \begin{enumerate}
  \item[2.] a concentric pair and a third circle exterior to the
  larger circle of the pair;
  \item[3.] one large circle containing two mutually exterior smaller circles,
  with the large circle concentric with one of them;
  \item[4.] three strictly nested circles with exactly two coincident centers:
  $4a$ for the outer and middle circles, $4b$ for the middle and inner
  circles, and $4c$ for the outer and inner circles.
  \end{enumerate}
\item[(III)] \emph{Three distinct collinear centers.}
  \begin{enumerate}
  \item[5.] three pairwise exterior circles;
  \item[6.] one contained pair and a third circle exterior to the outer circle;
  \item[7.] one large circle containing two mutually exterior small circles;
  \item[8.] three strictly nested circles.
  \end{enumerate}
\item[(IV)] \emph{Three non-collinear centers.}
  \begin{enumerate}
  \item[9.] three pairwise exterior circles;
  \item[10.] one contained pair and a third circle exterior to the outer circle;
  \item[11.] one large circle containing two mutually exterior small circles;
  \item[12.] three strictly nested circles.
  \end{enumerate}
\end{enumerate}
These four classes exhaust the possible arrangements of the
centers of three pairwise disjoint invariant circles. For three distinct
collinear centers, the circles are called coaxial if their pairwise radical
axes, namely the lines $f_i-f_j=0$, coincide, and non-coaxial otherwise.
The normal forms below cover all admissible containment configurations,
with their geometric parameter domains and circle labellings specified
in each case.

\begin{maintheorem}\label{main:concentric}
For configuration $1$ translate the common center to the origin.
The family of polynomial differential systems
of degree at most four having the three prescribed circles as invariant
algebraic curves is
\begin{equation}\label{sys:main-concentric}
 \dot x=-yq(x,y),\qquad \dot y=xq(x,y),
\end{equation}
where $q\in\R[x,y]$ is an arbitrary polynomial of total
degree at most three.
The polynomial $H=x^2+y^2$ is a global first integral
of system~\eqref{sys:main-concentric}. Every prescribed
circle either contains at least one equilibrium or belongs to a period annulus, and
hence none is an algebraic limit cycle.  Exact degree four requires
$\deg q=3$, whereas coprimality forces $q$ to be a non-zero constant.
\end{maintheorem}

\begin{maintheorem}\label{main:two-centers}
For configurations $2$--$4$ label the two concentric circles as
$f_1=0$ and $f_2=0$, and the remaining circle as $f_3=0$.
After a translation, a rotation and a homothety, write
\[
 f_1=x^2+y^2-r_1^2,\qquad f_2=x^2+y^2-r_2^2,\qquad
 f_3=(x-1)^2+y^2-r_3^2,
\]
where $r_i>0$ for $i=1,2,3$, $r_1\ne r_2$, and the three circles
are pairwise disjoint.  The family of real polynomial differential
systems of degree at most four having these circles as invariant algebraic
curves is
\begin{equation}\label{sys:main-two-centers}
 \dot x=yA,\qquad \dot y=\lambda f_1f_2-xA,
\end{equation}
where
\[
 A=\lambda f_1f_2-f_3B,\qquad
 B=\lambda(x^2+y^2)+\alpha x+\beta y+\gamma,
\]
with $\lambda,\alpha,\beta,\gamma\in\mathbb R$.
All three circles are algebraic limit cycles if and only if
\[
 \lambda\beta\ne0,\qquad
 |\lambda r_i^2+\gamma|>r_i\sqrt{\alpha^2+\beta^2},\qquad i=1,2.
\]
Whenever the above limit-cycle conditions hold,
system~\eqref{sys:main-two-centers} is coprime and of degree four,
and all three limit cycles are hyperbolic.
For $r_1>r_2>0$ and $r_3>0$, the strict geometric parameter domains are given
in Table~\ref{tab:two-centers-domains} in Section~\ref{sec:two-centers}.
Each domain admits a realization satisfying the above limit-cycle conditions.
\end{maintheorem}

\begin{maintheorem}\label{main:collinear}
For configurations $5$--$8$ normalize the three circles as
$f_i=0$, $i=1,2,3$, where
\[
 f_1=x^2+y^2-r_1^2,\qquad
 f_2=(x-1)^2+y^2-r_2^2,\qquad
 f_3=(x-u)^2+y^2-r_3^2,
\]
with $r_i>0$, $u\ne0,1$, and the three circles pairwise disjoint. Put
\[
 \Delta=u^2-u-r_3^2+ur_2^2+(1-u)r_1^2.
\]
The real polynomial differential systems of degree at most four having
these circles as invariant algebraic curves are as follows.
\begin{enumerate}
\item[\rm (i)] If $\Delta\ne0$ the family is
\begin{equation}\label{sys:main-collinear-noncoaxial}
\begin{aligned}
 \dot x&=y\bigl[(\lambda+\ell)f_1f_2
              -[\lambda+(1-u)\ell]f_2f_3-u\ell f_1f_3\bigr],\\
 \dot y&=[\lambda+(1-u)\ell]xf_2f_3
         +u\ell(x-1)f_1f_3-(\lambda+\ell)(x-u)f_1f_2,
\end{aligned}
\end{equation}
where $\ell=\alpha x+\beta y+\gamma$ and
$\lambda,\alpha,\beta,\gamma\in\mathbb R$ are arbitrary.
All three circles are algebraic limit cycles if and only if
$\lambda\beta\ne0$ and
\[
 |\lambda+(1-u)\gamma|>|1-u|r_1\sqrt{\alpha^2+\beta^2},\ 
 |\alpha+\gamma|>r_2\sqrt{\alpha^2+\beta^2},\ 
 |\lambda+\alpha u+\gamma|>r_3\sqrt{\alpha^2+\beta^2}.
\]
Whenever these conditions hold system~\eqref{sys:main-collinear-noncoaxial} is coprime and of exact degree
four, and all three limit cycles are hyperbolic.

\item[\rm (ii)] If $\Delta=0$ the family is
\begin{equation}\label{sys:main-collinear-coaxial}
 \dot x=q_2y(f_1-f_2),\qquad
 \dot y=q_2\bigl[(1-x)f_1+xf_2\bigr],
\end{equation}
where $q_2\in\mathbb R[x,y]$ is arbitrary with $\deg q_2\le2$.
The rational function $f_1/f_2$ is a first integral
of system~\eqref{sys:main-collinear-coaxial} wherever it is defined.
Every prescribed circle either contains an equilibrium or belongs to a
period annulus; hence none is a limit cycle. Exact degree four requires
$\deg q_2=2$, whereas coprimality forces $q_2$ to be a non-zero constant.
\end{enumerate}
Each of configurations $5$--$8$ admits a realization satisfying the
limit-cycle conditions in \textup{(i)}.
\par
\end{maintheorem}

For configurations $9$--$12$ normalize the three circles as
$f_i=0$, $i=1,2,3$, where
\[
 f_1=x^2+y^2-r_1^2,\qquad
 f_2=(x-1)^2+y^2-r_2^2,\qquad
 f_3=(x-u)^2+(y-v)^2-r_3^2,
\]
with $r_i>0$, $v>0$, and the three circles pairwise disjoint.
For real parameters $\lambda,\mu,\kappa$ introduce the following
quantities for the periodic-orbit and hyperbolicity criteria in the next
Theorem~\ref{main:noncollinear}:
\[
 d_{12}=1,\qquad d_{13}=\sqrt{u^2+v^2},\qquad
 d_{23}=\sqrt{(u-1)^2+v^2},\qquad d_{ik}=d_{ki}.
\]
Set
\[
 \omega_1=-\lambda-\mu-\kappa v,\qquad
 \omega_2=\lambda+\kappa v,\qquad
 \omega_3=\mu+\kappa v.
\]
For $k\ne i$ put $c_{ki}=d_{ki}^2+r_k^2-r_i^2$.
In the following definitions $\{i,j,k\}=\{1,2,3\}$.
When $|\omega_k|>|\kappa|r_kd_{ij}$ for $k=1,2,3$ define
\[
 \delta_k=\operatorname{sgn}(\omega_k)
 \sqrt{\omega_k^2-\kappa^2r_k^2d_{ij}^2},\qquad
 h_{ki}=\operatorname{sgn}(c_{ki})
             \sqrt{c_{ki}^2-4r_k^2d_{ki}^2},
\]
where the square roots are positive.
Set $H_{ki}=\omega_kc_{ki}+\delta_kh_{ki}-2\kappa vr_k^2$.
Define
\[
\begin{aligned}
 B_{12}&=-2\lambda-2u\mu-\kappa v(2u+c_{12}),\\
 B_{13}&=2u\lambda+2d_{13}^2\mu+\kappa v(2u+c_{13}),
\end{aligned}
\]
and let $B_{23}=-B_{12}-B_{13}$. Finally define
\[
\begin{aligned}
 \Phi_1&=B_{12}H_{13}+B_{13}H_{12},\\
 \Phi_2&=B_{12}H_{23}+B_{23}H_{21},\\
 \Phi_3&=B_{13}H_{32}+B_{23}H_{31}.
\end{aligned}
\]

\begin{maintheorem}\label{main:noncollinear}
All real polynomial differential systems of degree at most four
preserving the normalized circles above defined are
\begin{equation}\label{sys:main-noncollinear}
\begin{aligned}
 \dot x&=-yA_1f_2f_3-yA_2f_1f_3-(y-v)A_3f_1f_2,\\
 \dot y&=xA_1f_2f_3+(x-1)A_2f_1f_3+(x-u)A_3f_1f_2,
\end{aligned}
\end{equation}
where
\[
 A_1=-\lambda-\mu+\kappa[-vx+(u-1)y-v],\quad
 A_2=\lambda+\kappa(vx-uy),\quad
 A_3=\mu+\kappa y,
\]
and $\lambda,\mu,\kappa\in\mathbb R$ are arbitrary.
If $(\lambda,\mu,\kappa)\ne(0,0,0)$, system~\eqref{sys:main-noncollinear} has degree four.
\begin{enumerate}
\item[\rm (i)] All three circles are periodic orbits if and only if
\[
 |\lambda+\mu+\kappa v|>|\kappa|r_1d_{23},\quad
 |\lambda+\kappa v|>|\kappa|r_2d_{13},\quad
 |\mu+\kappa v|>|\kappa|r_3.
\]
Under these inequalities system~\eqref{sys:main-noncollinear} is coprime.
All three circles are hyperbolic algebraic limit cycles if and only if
\[
 \kappa\Phi_1\Phi_2\Phi_3\ne0.
\]

\item[\rm (ii)] If $\kappa=0$ every prescribed circle either consists entirely of equilibria
or belongs to a period annulus; hence none is a limit cycle.
\end{enumerate}
Each of configurations $9$--$12$ admits a realization with three
hyperbolic algebraic limit cycles.
\end{maintheorem}

\begin{re}\label{rem:noncollinear-nonhyperbolicity}
Assume the inequalities in \textup{(i)} and $\kappa\ne0$.
Then $\Phi_k=0$ if and only if the $k$-th circle is non-hyperbolic.
\end{re}

If a coprime system of exact degree four has the three prescribed circles
as limit cycles it has no fourth distinct invariant circle of positive
radius. Otherwise by \cite[Theorem~2]{LlibreRamirez2016} none of these
four invariant circles could be a limit cycle, contradicting the hypothesis.

An outline of this paper is as follows. Section~\ref{sec:preliminary}
collects preliminary results on circle geometry, characteristic exponents
and non-isolation of periodic orbits.
Sections~\ref{sec:concentric}--\ref{sec:noncollinear} treat the four classes
of circle centers in order and prove
Theorems~\ref{main:concentric}--\ref{main:noncollinear}, respectively.

\section{Preliminary results}
\label{sec:preliminary}

\subsection{Circle geometry and algebraic restrictions}
\label{subsec:circle-preliminaries}

For each $i=1,2,3$, let
$f_i(x,y)=(x-a_i)^2+(y-b_i)^2-r_i^2$, where $(a_i,b_i)$ and
$r_i>0$ are the center and radius of the $i$-th circle.
The value $f_i(x,y)$ is the power of the point $(x,y)$ with respect
to this circle. For two circles with distinct centers, their radical axis
is the line $f_i-f_j=0$, consisting of points having equal power with
respect to the two circles. This line is perpendicular to the line joining
their centers. A concentric pair of distinct circles has
no finite radical axis, since $f_i-f_j$ is a non-zero constant.

The complexification of a real circle is the affine curve defined by the
same equation in $\mathbb C^2$. Two such curves
$f_i=0$ and $f_j=0$ meet transversely at
$p$ if $\det(\nabla f_i(p),\nabla f_j(p))\ne0$, where
$\nabla f_i$ denotes the gradient of $f_i$.
By the complex inverse function theorem, we have local
holomorphic coordinates $(f_i,f_j)$ near $p$.

We record the algebraic facts used below; see also
\cite[Chapters~4--5]{CoxLittleOShea2015} for the relation between polynomial
ideals, affine varieties and rational functions.

\begin{lm}\label{lem:circle-restriction}
Let $f=(x-a)^2+(y-b)^2-r^2$, with $a,b\in\mathbb R$ and $r>0$.
The complex curve $f=0$ admits the parametrization
\begin{equation}\label{eq:complex-circle-parametrization}
 x(t)=a+\frac{t+r^2/t}{2},\qquad
 y(t)=b+\frac{t-r^2/t}{2\mathrm i},\qquad
 t\in\mathbb C\setminus\{0\}.
\end{equation}
The polynomial $f$ is irreducible in $\mathbb C[x,y]$.
For every polynomial $h\in\mathbb C[x,y]$ of total degree at most $m$,
where $m$ is a nonnegative integer, substitution of this parametrization
gives a Laurent polynomial of the form
\[
 h(x(t),y(t))=\sum_{k=-m}^{m}c_k t^k,\qquad c_k\in\mathbb C.
\]
A polynomial vanishes on the complex circle, or on its entire real circle,
if and only if it is divisible by $f$.
If $f_1,\ldots,f_N$ are distinct positive-radius circle polynomials then
they are pairwise relatively prime, and any polynomial vanishing on each
of their circles is divisible by $f_1\cdots f_N$.
\end{lm}

\begin{proof}
Set $u=x-a+\mathrm i(y-b)$, $v=x-a-\mathrm i(y-b)$.
This is an invertible affine change of variables, and $f=uv-r^2$.
On the curve $f=0$, set $t=u$. Since $uv=r^2\ne0$, we have
$t\ne0$ and $v=r^2/t$. Solving for $x,y$ gives
\eqref{eq:complex-circle-parametrization}. Conversely, these formulas
give a point on $f=0$ for every $t\ne0$.
The real circle is obtained by taking $t=re^{\mathrm i\theta}$,
$\theta\in\mathbb R$.

Since $uv-r^2$ has degree one in $v$, any nontrivial factorization
in $\mathbb C[u,v]$ would have a nonconstant factor depending only on
$u$. Such a factor would divide both coefficients $u$ and $-r^2$,
which is impossible because $r>0$. Thus $f$ is irreducible in
$\mathbb C[x,y]$.

Both $x(t)$ and $y(t)$ are linear combinations of $1,t,t^{-1}$.
Consequently, substituting them into a monomial $x^\alpha y^\beta$
gives only powers $t^k$ with $-(\alpha+\beta)\le k\le\alpha+\beta$.
Every monomial of $h$ has $\alpha+\beta\le m$, which proves the
stated Laurent polynomial representation.

Suppose that $h$ vanishes on the entire real circle. Polynomial
division with respect to $y-b$ gives
\[
 h=qf+A(x)(y-b)+B(x),
\]
where $q\in\mathbb C[x,y]$ and $A,B\in\mathbb C[x]$.
For each $x\in(a-r,a+r)$, evaluating this identity at the two points with
\[
 y-b=\pm\sqrt{r^2-(x-a)^2}
\]
gives $A(x)=B(x)=0$. Hence $A$ and $B$ are identically zero,
and $f$ divides $h$. Vanishing on the complex circle implies
vanishing on the real circle, while every multiple of $f$ vanishes
on the complex circle. This proves both equivalences.

Distinct circle polynomials are
non-associate irreducibles, and the product assertion follows from unique
factorization in $\mathbb C[x,y]$.
\end{proof}

\begin{lm}\label{lem:circle-complex-intersections}
Two disjoint real circles of positive radii with distinct centers
have exactly two distinct affine complex intersection points. Both are
non-real and transverse.
\end{lm}

\begin{proof}
A real translation and rotation put the centers at $(0,0)$ and
$(d,0)$, with $d>0$, without changing the complex intersections or
their transversality.
If the radii are $r,s>0$, the circle equations become
\[
 x^2+y^2=r^2,\qquad (x-d)^2+y^2=s^2.
\]
Subtracting gives $2dx=d^2+r^2-s^2$. Solving for $x$ and substituting
into the first equation yields
\[
 \begin{aligned}
 x&=\frac{d^2+r^2-s^2}{2d},\\
 y^2&=r^2-\left(\frac{d^2+r^2-s^2}{2d}\right)^2
     =-\frac{\bigl(d^2-(r+s)^2\bigr)\bigl(d^2-(r-s)^2\bigr)}{4d^2}.
 \end{aligned}
\]
Since the real circles are disjoint, either $d>r+s$ or
$0<d<|r-s|$. In either case the product in the numerator is positive,
so $y^2<0$. Thus $x$ is real and $y$ has exactly two distinct
non-zero purely imaginary values. These give exactly two distinct
non-real affine intersection points.
At either point the gradients of the two circle polynomials are
$(2x,2y)$ and $(2(x-d),2y)$, whose determinant is
$4dy\ne0$. Hence both intersections are transverse.
\end{proof}

\subsection{Characteristic exponents and cofactors}
\label{subsec:characteristic-exponents}

Let $\Gamma$ be a periodic orbit of a $C^1$ planar vector field
$\vf=(P,Q)$, parametrized by a solution $\gamma(t)$ of minimal period
$T>0$.  Let $\Pi$ be its positive-time Poincar\'e return map in a local
transverse coordinate with $\Pi(0)=0$.  The derivative formula
\cite[Theorem~1.23]{DumortierLlibreArtes2006} gives
\begin{equation}\label{eq:preliminary-multiplier}
 \Pi'(0)=\exp\!\left(\int_0^T
 \operatorname{div}\vf(\gamma(t))\,dt\right),\qquad
 \operatorname{div}\vf=P_x+Q_y.
\end{equation}
Throughout this paper the \emph{characteristic exponent} means the
unnormalized logarithm
$\chi(\Gamma)=\log\Pi'(0)=\int_0^T\operatorname{div}\vf(\gamma(t))\,dt$.
Thus $\Pi'(0)=e^{\chi(\Gamma)}$ is the characteristic multiplier.

Hyperbolicity is defined for periodic orbits, whereas a
limit cycle is an isolated periodic orbit.
Consequently $\Gamma$ is hyperbolic if and only if
$\chi(\Gamma)\ne0$, equivalently $\Pi'(0)\ne1$;
see \cite[Section~2.2, p.~48]{Zhang2017}.
In that case it is an isolated periodic orbit,
and hence a hyperbolic limit cycle,
attracting when $\chi(\Gamma)<0$ and repelling when
$\chi(\Gamma)>0$; see \cite[Theorem~1]{GiacominiGrau2005}.
Indeed $(\Pi-\mathrm{id})'(0)=e^{\chi(\Gamma)}-1\ne0$, so the fixed
point is isolated. When $\chi(\Gamma)=0$ the orbit is
non-hyperbolic: it may be a non-hyperbolic limit cycle or a non-isolated
periodic orbit, for example one belonging to a period annulus.
The equality alone does not determine its isolation or stability.

Suppose also that $\Gamma\subset\{f=0\}$, where $f$ is $C^1$,
$\nabla f\ne0$ on $\Gamma$, and $\vf(f)=Kf$ near $\Gamma$
with a $C^1$ cofactor $K$.  Then
\cite[Theorem~2]{GiacominiGrau2005} yields
\begin{equation}\label{eq:preliminary-cofactor}
 \chi(\Gamma)=\int_0^T\operatorname{div}\vf(\gamma(t))\,dt
             =\int_0^T K(\gamma(t))\,dt.
\end{equation}
The regularity hypothesis holds for every circle of positive radius.
The equality concerns the integrals over a full period and does not assert
that $K=\operatorname{div}\vf$ pointwise.

For a periodic circle parametrized by $\theta\in[0,2\pi]$ write its
angular equation as $\dot\theta=W(\theta)$.
Combining the multiplier formula
\cite[Theorem~1.23]{DumortierLlibreArtes2006} with the cofactor identity
\cite[Theorem~2]{GiacominiGrau2005} and using $dt=d\theta/W(\theta)$,
with $W$ of constant non-zero sign and the positive-time orientation
taken into account, gives
\begin{equation}\label{eq:preliminary-angle-exponent}
 I=\int_0^{2\pi}\frac{K(\theta)}{W(\theta)}\,d\theta,\qquad
 \chi(\Gamma)=\operatorname{sgn}(W)I,
\end{equation}
where $K(\theta)$ denotes the cofactor restricted to the circle.
The sign corrects for clockwise motion when $W<0$; all stability
statements use the positive-time exponent $\chi(\Gamma)$.

\subsection{Criteria for non-isolation}
\label{subsec:non-isolation}

The following two criteria will be used to exclude isolated periodic
orbits. For the general theory of planar flows and first integrals see
\cite{DumortierLlibreArtes2006}.

\begin{lm}\label{lem:regular-first-integral}
Let $\Gamma$ be a periodic orbit of a $C^1$ planar vector field.
Suppose a single-valued $C^1$ first integral $H$ is defined in a
neighborhood of $\Gamma$ and satisfies $\nabla H\ne0$ on $\Gamma$.
Then $\Gamma$ belongs to a local period annulus and is not a limit cycle.
\end{lm}

\begin{proof}
Let $h_0$ be the constant value of $H$ on $\Gamma$.
By continuity and compactness, choose an annular neighborhood $U$
of $\Gamma$ on which both $\nabla H$ and the vector field are non-zero.
If $\Gamma$ were a limit cycle, trajectories in a sufficiently small
annular region on one side of $\Gamma$ would spiral toward $\Gamma$
in forward or backward time. Constancy of $H$ along trajectories
and continuity at $\Gamma$ would then give $H\equiv h_0$ throughout
that region, contradicting $\nabla H\ne0$ on $U$.

By the implicit function theorem and compactness of $\Gamma$,
a smaller annular neighborhood is foliated by nearby simple closed
level curves of $H$. The vector field is tangent to these curves
and has no zeros there, so each curve is a periodic orbit.
Hence $\Gamma$ belongs to a local period annulus.
\end{proof}

For the reflection $R(x,y)=(x,-y)$ a vector field is reversible if
$\vf(Rz)=-DR\,\vf(z)$, where $z=(x,y)^T$; equivalently $P$
is odd and $Q$ is even in $y$. Its flow satisfies
$R\circ\phi_t=\phi_{-t}\circ R$.

\begin{lm}\label{lem:reversible-period-annulus}
Let a $C^1$ planar vector field be reversible with respect to reflection
in a line. A periodic orbit meeting that line transversely at two distinct
points belongs to a local period annulus and is not a limit cycle.
\end{lm}

\begin{proof}
By a rigid change of coordinates take the line to be $y=0$.
Transversality and continuous dependence on initial conditions imply that
trajectories starting near the first intersection on this line meet it
again near the second intersection. Reversibility reflects each such arc
into the remainder of a periodic orbit, by uniqueness of solutions.
These nearby periodic orbits form a local annulus. This is the planar
application of the symmetric-orbit criterion in
\cite[Theorem~4.1]{LambRoberts1998}.
\end{proof}

Multiplying a vector field by a real function which is non-zero in a
neighborhood preserves its trajectories there, possibly reversing their
orientation. Hence it also preserves the existence of a local period annulus.
\par

\section{Three coincident centers: configuration 1}
\label{sec:concentric}

This section proves Theorem~\ref{main:concentric} by first establishing a
criterion for concentric invariant circles valid for arbitrary degree.
The cubic case with two concentric invariant circles is given in
\cite[Theorem~1.1]{GineLlibreValls2018}; a related result when the degree
equals the number of circles appears in
\cite[Proposition~8]{LlibreRamirez2016}. Bounds for the number of concentric
invariant circles, counted with multiplicity and under the assumption that
there are only finitely many invariant circles, were obtained in
\cite[Theorem~1.1(b) and Corollary~1.1(b)]{BolanosLlibre2011}.

\begin{prop}\label{prop:concentric-higher-degree}
Let $n,m\ge1$ be integers with $2m>n$, and let
\[
 f_i=x^2+y^2-r_i^2,\qquad i=1,\ldots,m,\qquad
 0<r_1<\cdots<r_m.
\]
The complete family of real polynomial systems of degree at most $n$
having the circles $f_i=0$ invariant is
\[
 \dot x=-yq(x,y),\qquad \dot y=xq(x,y),
\]
where $q\in\mathbb R[x,y]$ is an arbitrary polynomial of total degree at
most $n-1$. The polynomial $H=x^2+y^2$ is a global first integral,
and all prescribed cofactors vanish. Every circle of positive radius
centered at the origin is invariant and either contains an equilibrium
or belongs to a period annulus; hence none is an algebraic limit cycle.
\end{prop}

\begin{proof}
Let $G=\prod_{i=1}^m f_i$. The positive, distinct radii make the circle
polynomials pairwise relatively prime. Invariance gives
$f_i\mid xP+yQ$ for every $i$, so $G\mid xP+yQ$.
Since $2m>n$ we have
\[
 \deg(xP+yQ)\le n+1\le2m=\deg G,
\]
and therefore $xP+yQ=cG$ for a real constant $c$. At the origin
\[
 0=cG(0,0)=c(-1)^m\prod_{i=1}^m r_i^2,
\]
which forces $c=0$. As $\gcd(x,y)=1$, the identity
$xP+yQ=0$ gives $P=-yq$, $Q=xq$, with $\deg q\le n-1$.
Conversely, every such $q$ yields $\dot H=0$, so every circle centered
at the origin is invariant and all prescribed cofactors are zero.

If $q$ vanishes at a point of a circle of positive radius, that point
is an equilibrium. Otherwise, continuity and compactness give a surrounding
annulus on which $q$ is nowhere zero. The angular velocity on each
concentric circle in this annulus is $\dot\theta=q$, so these circles
are periodic orbits and form a period annulus.
\end{proof}

For the quartic case, translate the common center to the origin and write
the three prescribed circles as
\begin{equation}\label{eq:three-circles}
 f_i=x^2+y^2-r_i^2=0,\qquad i=1,2,3,
\end{equation}
where $0<r_1<r_2<r_3$.
Taking $n=4$ and $m=3$ in Proposition~\ref{prop:concentric-higher-degree}
gives the following corollary; the degree and coprimality assertions use
$\deg(-yq,xq)=\deg q+1$ and $\gcd(-yq,xq)\doteq q$ for
$q\not\equiv0$.

\begin{cor}\label{cor:concentric-quartic}
A real polynomial system of degree at most four has the three circles
\eqref{eq:three-circles} as invariant algebraic curves if and only if
\begin{equation}\label{sys:normal-form}
 \dot x=-yq(x,y),\qquad \dot y=xq(x,y),
\end{equation}
where $q\in\mathbb R[x,y]$ is an arbitrary polynomial of total degree at
most three. The polynomial $H=x^2+y^2$ is a global first integral, and
the three prescribed cofactors are zero. Every circle of positive radius
centered at the origin is invariant and either contains an equilibrium
or belongs to a period annulus. In particular, none of the three prescribed
circles is a limit cycle.
The system has exact degree four if and only if $\deg q=3$.
If $q\not\equiv0$, then $\gcd(P,Q)\doteq q$. Consequently, the system
is coprime if and only if $q$ is a non-zero constant, in which case it is
a linear center. In particular, no coprime system of exact degree four
has three distinct concentric invariant circles.
\end{cor}

For this family, the equilibria are the origin and the zeros of $q$, and
every non-equilibrium orbit lies on a circle centered at the origin.
Corollary~\ref{cor:concentric-quartic} proves Theorem~\ref{main:concentric}.

\section{Exactly two coincident centers: configurations 2, 3 and 4}
\label{sec:two-centers}

This section proves Theorem~\ref{main:two-centers}.

We now assume that exactly two circle centers coincide.  A translation, a
rotation and a homothety put the common center at the origin and the
remaining center at $(1,0)$.  Label the two concentric circles as $f_1=0$ and $f_2=0$,
and the remaining circle as $f_3=0$.  Their equations are
\begin{equation}\label{eq:two-centers-circles}
 f_1=x^2+y^2-r_1^2,\qquad f_2=x^2+y^2-r_2^2,\qquad f_3=(x-1)^2+y^2-r_3^2,
\end{equation}
where $r_i>0$ for $i=1,2,3$, $r_1\ne r_2$,
and $f_3=0$ has no common point with either $f_1=0$ or $f_2=0$.  Equivalently
for each $i=1,2$ one and only one of the two inequalities
\begin{equation}\label{eq:two-centers-disjointness}
 1>r_i+r_3,\qquad 1<|r_i-r_3|
\end{equation}
holds.  The inequalities are strict, so intersection and tangency are both
excluded.

For the geometric classification the concentric pair of circles we can assume
$r_1>r_2$. Table~\ref{tab:two-centers-domains} gives the strict domains;
in configuration 4 the three rows distinguish the position of $D_3$
in the nested chain.

\begin{table}[htbp]
\centering
\renewcommand{\arraystretch}{1.3}
\begin{tabular}{c|l|l}
Configuration & Geometry & Strict inequalities \\ \hline
2 & $D_2\subset\operatorname{int}D_1$, $D_3$ exterior to $D_1$
  & $r_1+r_3<1$ \\
3 & $D_2,D_3\subset\operatorname{int}D_1$, mutually exterior
  & $1+r_3<r_1,\ r_2+r_3<1$ \\
4a & $D_3\subset\operatorname{int}D_2\subset\operatorname{int}D_1$
  & $1+r_3<r_2$ \\
4b & $D_2\subset\operatorname{int}D_1\subset\operatorname{int}D_3$
  & $1+r_1<r_3$ \\
4c & $D_2\subset\operatorname{int}D_3\subset\operatorname{int}D_1$
  & $1+r_2<r_3,\ 1+r_3<r_1$ \\
\end{tabular}
\par\medskip
\caption{Strict domains for the configurations 2--4, with $r_i>0$ for
$i=1,2,3$, $r_1>r_2$, and $D_1,D_2$ concentric.}
\label{tab:two-centers-domains}
\end{table}

For the real parameters $\lambda,\alpha,\beta,\gamma$, we define
\begin{equation}\label{eq:B-two-centers}
 B(x,y)=\lambda (x^2+y^2)+\alpha x+\beta y+\gamma
\end{equation}
and
\begin{equation}\label{eq:A-two-centers}
 A(x,y)=\lambda f_1 f_2-f_3B.
\end{equation}

\begin{prop}\label{prop:two-centers-normal-form}
Assume that $r_i>0$ for $i=1,2,3$, $r_1\ne r_2$, and
\eqref{eq:two-centers-disjointness} holds.  A real polynomial differential system of
degree at most four has the three circles \eqref{eq:two-centers-circles} as
invariant algebraic curves if and only if it is
\begin{equation}\label{sys:two-centers-normal-form}
 \dot x=yA,\qquad \dot y=\lambda f_1f_2-xA.
\end{equation}
The corresponding cofactors are
\begin{equation}\label{eq:two-centers-cofactors}
 K_1=2\lambda yf_2,\qquad
 K_2=2\lambda yf_1,\qquad
 K_3=2yB.
\end{equation}
\end{prop}

\begin{proof}
For the two concentric circles we have
\[
 \vf(f_1)=\vf(f_2)=2(xP+yQ).
\]
Their invariance and coprimality imply that $f_1f_2$ divides
$\vf(f_1)$.  Since $\deg\vf(f_1)\leq5$ we can write
\begin{equation}\label{eq:two-centers-derivative-affine}
 \vf(f_1)=\vf(f_2)=2f_1f_2(ux+vy+w)
\end{equation}
for some constants $u,v,w$.  Invariance of the third circle gives
\[
 \vf(f_3)=\vf(f_1)-2P=K_3f_3,\qquad \deg K_3\leq3.
\]
Hence $P=f_1f_2(ux+vy+w)-\tfrac12f_3K_3$.  Substituting this in
$\vf(f_1)=2(xP+yQ)$ and setting $y=0$ yields
\begin{equation}\label{eq:axis-divisibility}
 xf_3(x,0)K_3(x,0)=2(x-1)f_1(x,0)f_2(x,0)(ux+w).
\end{equation}
At $x=0$ this identity and $r_1r_2>0$ imply $w=0$.

The roots of $f_3(x,0)$ are $1\pm r_3$.
Strict disjointness implies that neither is a root of
$f_1(x,0)f_2(x,0)$.  In particular at $x=1+r_3>0$
all factors on the right of \eqref{eq:axis-divisibility} except
$u$ are nonzero, so $u=0$.  Renaming $v=\lambda$ we obtain
\begin{equation}\label{eq:two-centers-derivative-final}
 \vf(f_1)=\vf(f_2)=2\lambda yf_1f_2.
\end{equation}
Equation \eqref{eq:axis-divisibility} now gives $K_3(x,0)=0$, so
$K_3=2yB$ for a polynomial $B$ of degree at most two.  Thus
\[
 P=y(\lambda f_1f_2-f_3B).
\]
Since $\deg P\leq4$ the homogeneous degree-four part inside
parentheses must vanish.  If $B_2$ is the homogeneous quadratic
part of $B$ this says
\[
 \lambda(x^2+y^2)^2-(x^2+y^2)B_2=0,
\]
so $B_2=\lambda(x^2+y^2)$.  Therefore $B$ has the form
\eqref{eq:B-two-centers}.  Recovering $Q$ from
$\vf(f_1)=2(xP+yQ)=2\lambda yf_1f_2$ gives
\eqref{sys:two-centers-normal-form}.

Conversely the displayed system satisfies
\[
 \vf(f_1)=2\lambda yf_1f_2,\qquad
 \vf(f_2)=2\lambda yf_1f_2,\qquad
 \vf(f_3)=2yf_3B.
\]
These identities prove invariance and the cofactor formulas
\eqref{eq:two-centers-cofactors}.
\end{proof}

\begin{prop}\label{prop:two-centers-degree-gcd}
For the family \eqref{sys:two-centers-normal-form}:
\begin{enumerate}
\item the system has exact degree four if and only if
\begin{equation}\label{eq:two-centers-exact-degree}
 (\lambda,\alpha,\beta)\ne(0,0,0);
\end{equation}
\item the components $P,Q$ are relatively prime if and only if
\begin{equation}\label{eq:two-centers-coprime}
 \lambda\ne0,\qquad
 (\alpha,\beta,\gamma+\lambda r_1^2)\ne(0,0,0),\qquad
 (\alpha,\beta,\gamma+\lambda r_2^2)\ne(0,0,0).
\end{equation}
\end{enumerate}
\end{prop}

\begin{proof}
The homogeneous quartic parts of \eqref{sys:two-centers-normal-form} are
\begin{equation}\label{eq:two-centers-quartic-parts}
 \begin{split}
 P_4&=y(x^2+y^2)\big((2\lambda-\alpha)x-\beta y\big),\\
 Q_4&=(x^2+y^2)\big((\alpha-\lambda)x^2+\beta xy+\lambda y^2\big).
 \end{split}
\end{equation}
They vanish simultaneously if and only if
$\lambda=\alpha=\beta=0$, which proves the first statement.

For the second statement if $\lambda=0$ then $P=-yf_3B$ and
$Q=xf_3B$, so $f_3$ is a common factor.
Suppose $\lambda\ne0$.  Any common irreducible factor of $P,Q$ divides
$xP+yQ=\lambda yf_1f_2$.  The factor $y$ cannot divide $Q$ when
$\lambda\ne0$, and $f_3$ is relatively prime to $yf_1f_2$.  The factor $f_1$
divides both components exactly when $f_1\mid B$, that is, when
$B=\lambda f_1$; this is the first excluded equality in
\eqref{eq:two-centers-coprime}.  Similarly $f_2$ is common exactly when
$B=\lambda f_2$.  There are no remaining possible common factors.
\end{proof}

\subsection{Periodic-orbit and limit-cycle conditions}

For $i=1,2$ parametrize $f_i=0$ by
$x=r_i\cos\theta$, $y=r_i\sin\theta$.  On $f_3=0$ use
$x=1+r_3\cos\varphi$, $y=r_3\sin\varphi$.  The angular equations are
\begin{equation}\label{eq:two-centers-angular}
 \left.\dot\theta\right|_{f_1=0}=f_3B,\qquad
 \left.\dot\theta\right|_{f_2=0}=f_3B,\qquad
 \left.\dot\varphi\right|_{f_3=0}=-\lambda f_1f_2.
\end{equation}
Since the three circles are pairwise disjoint the restrictions of $f_1,f_2,f_3$
appearing on the right-hand sides never vanish.  Moreover for $i=1,2$
\begin{equation}\label{eq:B-circle-test}
 B(r_i\cos\theta,r_i\sin\theta)=\lambda r_i^2+\gamma
 +r_i(\alpha\cos\theta+\beta\sin\theta).
\end{equation}
Thus for each $i=1,2$ $B$ has no zero on $f_i=0$ if and only if
\begin{equation}\label{eq:B-circle-inequality}
 |\lambda r_i^2+\gamma|>r_i\sqrt{\alpha^2+\beta^2}.
\end{equation}

\begin{prop}\label{prop:two-centers-limit-cycles}
All three circles in \eqref{eq:two-centers-circles} are algebraic limit cycles
of system \eqref{sys:two-centers-normal-form} if and only if
\begin{equation}\label{eq:three-limit-cycle-conditions}
 \quad
 \lambda\beta\ne0,\qquad
 |\lambda r_1^2+\gamma|>r_1\sqrt{\alpha^2+\beta^2},\qquad
 |\lambda r_2^2+\gamma|>r_2\sqrt{\alpha^2+\beta^2}.
 \quad
\end{equation}
Whenever \eqref{eq:three-limit-cycle-conditions} holds the system is
automatically coprime and of exact degree four, and all three limit cycles are
hyperbolic.
\end{prop}

\begin{proof}
By \eqref{eq:two-centers-angular} the three invariant circles are
periodic orbits if and only if $\lambda\ne0$ and $B$ has no zero
on either $f_1=0$ or $f_2=0$.  The latter requirement is precisely
\eqref{eq:B-circle-inequality} for $i=1,2$.

It remains to determine isolation, using the cofactor formula
\eqref{eq:preliminary-cofactor} and the orientation convention
\eqref{eq:preliminary-angle-exponent}.
On $f_1=0$ the cofactor is
$K_1=2\lambda r_1(r_1^2-r_2^2)\sin\theta$.  Hence
\[
 I_1=2\lambda r_1(r_1^2-r_2^2)
 \int_0^{2\pi}\frac{\sin\theta\,d\theta}
 {f_{3,1}(\theta)\bigl(B_1(\theta)+\beta r_1\sin\theta\bigr)},
\]
where $f_{3,1}(\theta)=r_1^2+1-r_3^2-2r_1\cos\theta$ and
$B_1(\theta)=\lambda r_1^2+\gamma+\alpha r_1\cos\theta$.
Split the integral at $\pi$ and replace $\theta$ by
$2\pi-\theta$ in the second half.  Both $f_{3,1}$ and $B_1$
are unchanged by this replacement, whereas the sine changes sign.
After factoring out $1/f_{3,1}(\theta)$ we obtain the following
identity for the paired terms:
\[
 \frac{\sin\theta}{B_1(\theta)+\beta r_1\sin\theta}
 -\frac{\sin\theta}{B_1(\theta)-\beta r_1\sin\theta}
 =\frac{-2\beta r_1\sin^2\theta}
 {B_1(\theta)^2-\beta^2r_1^2\sin^2\theta}.
\]
This gives
\begin{equation}\label{eq:chi-1}
 I_1=-4\lambda(r_1^2-r_2^2)\beta r_1^2
 \int_0^\pi\frac{\sin^2\theta\,d\theta}
 {f_{3,1}(\theta)
  \left(B_1(\theta)^2-\beta^2r_1^2\sin^2\theta\right)}.
\end{equation}
The same formula with $r_1$ and $r_2$ interchanged gives
$I_2$.  Under the disjointness and periodicity assumptions
$f_{3,1}$ is sign-definite and
$B_1(\theta)^2-\beta^2r_1^2\sin^2\theta>0$, since the latter expression
is the product of two values of $B$ with the same non-zero sign
at the paired points of the circle.
Hence $I_1$ and $I_2$ are non-zero
exactly when $\beta\ne0$.  For the third circle
\begin{equation}\label{eq:chi-3}
 I_3=-\frac{2\beta}{\lambda}
 \int_0^{2\pi}\frac{r_3^2\sin^2\varphi}
 {f_{1,3}(\varphi)f_{2,3}(\varphi)}\,d\varphi,
\end{equation}
where $f_{j,3}(\varphi)=f_j(1+r_3\cos\varphi,r_3\sin\varphi)$
for $j=1,2$.  This formula follows because the part of $B$ which is even in $\sin\varphi$ integrates to
zero.  The product $f_{1,3}f_{2,3}$ is sign-definite, so $I_3\ne0$ exactly
when $\beta\ne0$.  Thus by Section~\ref{sec:preliminary}
$\beta\ne0$ makes all three periodic orbits hyperbolic and isolated.

Let $\chi_i$ denote the characteristic exponent of the
periodic circle $f_i=0$, for $i=1,2,3$.
By \eqref{eq:preliminary-angle-exponent}
$\chi_i=\operatorname{sgn}(W_i)I_i$, where $W_i$ is the angular
velocity in \eqref{eq:two-centers-angular}; hence $\chi_i=0$ if and
only if $I_i=0$.

If $\beta=0$ then $P$ is odd and $Q$ is even in
$y$, so the vector field is reversible with respect to $y=0$.
Each prescribed periodic circle meets this line transversely at two points;
Lemma~\ref{lem:reversible-period-annulus} therefore excludes isolation.
This proves necessity.
Finally $\lambda\beta\ne0$ implies both
\eqref{eq:two-centers-exact-degree} and \eqref{eq:two-centers-coprime}.
\end{proof}

\subsection{Hyperbolic realizations}

For any geometric parameters in Table~\ref{tab:two-centers-domains} use
the following parameter choice:
\begin{equation}\label{eq:universal-example-parameters}
 \lambda=1,\qquad \alpha=0,\qquad \beta=1,\qquad
 \gamma>\max_{1\leq i\leq3}\{r_i-r_i^2\}
\end{equation}
These parameters satisfy \eqref{eq:three-limit-cycle-conditions}, so
Proposition~\ref{prop:two-centers-limit-cycles} yields a coprime quartic
system with three hyperbolic circular limit cycles in each configuration.

This completes the proof of Theorem~\ref{main:two-centers}.

\section{Three distinct collinear centers: configurations 5, 6, 7 and 8}
\label{sec:collinear}

This section proves Theorem~\ref{main:collinear}.

Apply a translation, a rotation and a homothety and reflect
the line of centers if needed. Write
\begin{equation}\label{eq:collinear-circles}
 f_1=x^2+y^2-r_1^2,\qquad
 f_2=(x-1)^2+y^2-r_2^2,\qquad
 f_3=(x-u)^2+y^2-r_3^2,
\end{equation}
where $r_1,r_2,r_3>0$ and $u\ne0,1$.
The strict geometric domains in Table~\ref{tab:collinear-domains}
ensure that the three real circles are pairwise disjoint. The disk labels
follow the containment roles indicated in each row.

\begin{table}[htbp]
\centering
\renewcommand{\arraystretch}{1.3}
\begin{tabular}{c|l|l}
Configuration & Geometry & Parameter conditions \\ \hline
5 & Pairwise exterior disks
  & $u\ge2,\ r_1+r_2<1,$ \\
  & & $r_2+r_3<u-1$ \\
6 & $D_2\subset\operatorname{int}D_1$, $D_3$ exterior to $D_1$
  & $1+r_2<r_1,\ r_1+r_3<|u|$ \\
7 & $D_2,D_3\subset\operatorname{int}D_1$, mutually exterior
  & $|u|\ge1,\ 1+r_2<r_1,$ \\
  & & $|u|+r_3<r_1,\ r_2+r_3<|u-1|$ \\
8 & $D_3\subset\operatorname{int}D_2\subset\operatorname{int}D_1$
  & $1+r_2<r_1,\ |u-1|+r_3<r_2$ \\
\end{tabular}
\par\medskip
\caption{Domains for configurations 5--8, with $r_i>0$ for $i=1,2,3$,
$u\in\mathbb{R}\setminus\{0,1\}$, and strict disjointness inequalities.}
\label{tab:collinear-domains}
\end{table}

In configuration 5 order the centers along the line, normalize the smaller
adjacent-center gap to one, and use reflection to place it on the left;
when $u=2$ impose the residual convention $r_1\le r_3$.
In configuration 7 label the nearer small-circle center as $C_2$,
normalize its distance from $C_1$ to one, and reflect it to the positive
axis; when $u=-1$ impose $r_2\le r_3$.
In configurations 6 and 8 normalize the distance between $C_1$ and
$C_2$ to one and reflect the line so that $C_2=(1,0)$.
The remaining exterior or containment inequalities follow from those in
the table and the triangle inequality.

Introduce the quantity
\begin{equation}\label{eq:Delta-collinear}
 \Delta:=u^2-u-r_3^2+u r_2^2+(1-u)r_1^2.
\end{equation}
Then we have
\begin{equation}\label{eq:pencil-relation}
 (1-u)f_1+u f_2-f_3=-\Delta.
\end{equation}
The radical-axis definition in Section~\ref{subsec:circle-preliminaries}
gives the following equations for the axes of the pairs $(1,2)$ and
$(1,3)$:
\[
 x=\frac{1+r_1^2-r_2^2}{2},\qquad
 x=\frac{u^2+r_1^2-r_3^2}{2u},
\]
respectively. The difference between the two right-hand sides is
\[
 \frac{u^2+r_1^2-r_3^2}{2u}
 -\frac{1+r_1^2-r_2^2}{2}=\frac{\Delta}{2u}.
\]
Since $u\ne0$ these two axes coincide exactly when $\Delta=0$.
If they coincide then $f_2=f_3$ along their common line, so the
third radical axis coincides with them as well.

Thus $\Delta=0$ if and only if the three circles are coaxial.
The cases $\Delta=0$ and $\Delta\ne0$ lead to different
normal forms and are treated separately below.

Set
\begin{equation}\label{eq:T-collinear}
 \mathcal T_1=f_2f_3\binom{-y}{x},\qquad
 \mathcal T_2=f_1f_3\binom{-y}{x-1},\qquad
 \mathcal T_3=f_1f_2\binom{-y}{x-u},
\end{equation}
and define the polynomial vector fields
\begin{equation}\label{eq:Z-collinear}
 \mathcal Z_0=(1-u)\mathcal T_1+u\mathcal T_2-\mathcal T_3,
 \qquad
 \mathcal Z_1=\mathcal T_1-\mathcal T_3.
\end{equation}
Although every $\mathcal T_i$ has degree five the degree-five terms cancel
in both fields in \eqref{eq:Z-collinear}.
Moreover the degree-four terms also cancel in $\mathcal Z_0$.  Consequently we have $\deg\mathcal Z_0=3$ and
$\deg\mathcal Z_1=4$.
The fields $\mathcal Z_0$ and $\mathcal Z_1$ will be used
to express the non-coaxial normal form in Proposition~\ref{prop:collinear-normal-form}.

\subsection{A representation with affine coefficients}
\label{subsec:uniform-affine-representation}

For general methods of representing polynomial vector fields with prescribed
invariant algebraic curves see \cite{ChristopherLlibrePantaziWalcher2009}.

The following lemma applies both to the collinear,
non-coaxial case $\Delta\ne0$ considered here and to the non-collinear
case treated in Section~\ref{sec:noncollinear}. For
$C_1=(0,0)^T,C_2=(1,0)^T,C_3=(u,0)^T$ its fields $\mathcal T_i$
coincide with those in \eqref{eq:T-collinear}. The hypothesis that the
three complexified circles have no common finite point is verified after
the proof. The proof obtains the affine coefficients by
restricting the vector field to each complex circle and comparing the
powers of $t$ occurring in the resulting Laurent polynomials.

\begin{lm}\label{lem:uniform-affine-representation}
Let $f_i=(x-a_i)^2+(y-b_i)^2-r_i^2$, $i=1,2,3$, with $r_i>0$.
Assume that the circles $f_i=0$ are pairwise disjoint and have pairwise
distinct centers $C_i=(a_i,b_i)^T$.
Suppose the three complexified circles have no common finite complex point.
Set $z=(x,y)^T$ and let $J$ be the linear map representing
counterclockwise rotation through $90^\circ$, so that
$J(x,y)^T=(-y,x)^T$. Define
\[
 \mathcal T_i=\left(\prod_{j\ne i}f_j\right)J(z-C_i).
\]
Then the following statements hold.
\begin{enumerate}
\item[(i)] Every real polynomial vector field $\vf$ of degree at most four
leaving the three circles invariant can be written as
\begin{equation}\label{eq:uniform-affine-representation}
 \vf=\sum_{i=1}^3 A_i\mathcal T_i,
\end{equation}
with $A_i=\ell_i+c_i$, where the $\ell_i$ are real homogeneous linear
polynomials and the $c_i$ are
real constants.
\item[(ii)] The coefficients satisfy
\begin{equation}\label{eq:uniform-affine-constraints}
 \sum_{i=1}^3\ell_i=0,\qquad
 \sum_{i=1}^3 C_i(\nabla\ell_i)^T
 +\Bigl(\sum_{i=1}^3 c_i\Bigr)\mathrm{Id}_2=0.
\end{equation}
\item[(iii)] Conversely any such coefficients satisfying
\eqref{eq:uniform-affine-constraints} define, through
\eqref{eq:uniform-affine-representation}, a vector field of degree at most
four leaving all three circles invariant.
\item[(iv)] The coefficients $A_1,A_2,A_3$ in
\eqref{eq:uniform-affine-representation} are unique.
\end{enumerate}
\end{lm}

\begin{proof}
Lemma~\ref{lem:circle-complex-intersections} and the
absence of a common finite complex point give four distinct transverse affine
intersection points on each complexified circle.

Fix $i$ and denote by $(x_i(t),y_i(t))$ the parametrization
\eqref{eq:complex-circle-parametrization} with $(a,b,r)=(a_i,b_i,r_i)$.
Put $P_i(t)=P(x_i(t),y_i(t))$, $Q_i(t)=Q(x_i(t),y_i(t))$, and define
\[
 S_i(t)=\frac{(x_i(t)-a_i)Q_i(t)-(y_i(t)-b_i)P_i(t)}{r_i^2}.
\]
These are Laurent polynomials. Tangency, namely
$(x_i(t)-a_i)P_i(t)+(y_i(t)-b_i)Q_i(t)=0$, gives
$P_i(t)=-(y_i(t)-b_i)S_i(t)$ and
$Q_i(t)=(x_i(t)-a_i)S_i(t)$, hence
\[
 P_i(t)+\mathrm iQ_i(t)=\mathrm i tS_i(t),\qquad
 P_i(t)-\mathrm iQ_i(t)=-\mathrm i r_i^2t^{-1}S_i(t).
\]
By Lemma~\ref{lem:circle-restriction} and $\deg P,\deg Q\le4$,
both $P_i(t)$ and $Q_i(t)$ are linear combinations of
$t^{-4},\ldots,t^4$. If $t^\ell$ has a non-zero coefficient in
$S_i(t)$, the first identity gives $-4\le\ell+1\le4$, whereas
the second gives $-4\le\ell-1\le4$. Thus $-3\le\ell\le3$, and
\[
 S_i(t)=\sum_{\ell=-3}^{3}s_{i\ell}t^\ell,\qquad s_{i\ell}\in\mathbb C.
\]

At each intersection with another circle the two independent tangency
conditions force $\vf=0$, and thus $S_i=0$.
For $j\ne i$, set $G_{ij}(t)=f_j(x_i(t),y_i(t))$.
On the circle $f_i=0$ we have $f_j=f_j-f_i$, and this difference
is an affine polynomial. More explicitly,
\[
G_{ij}(t)=\alpha_{ij}t+\beta_{ij}+\gamma_{ij}t^{-1},
\]
where
\[
 \begin{aligned}
 \alpha_{ij}&=a_i-a_j-\mathrm i(b_i-b_j),\\
 \beta_{ij}&=(a_i-a_j)^2+(b_i-b_j)^2+r_i^2-r_j^2,\\
 \gamma_{ij}&=r_i^2\bigl(a_i-a_j+\mathrm i(b_i-b_j)\bigr).
 \end{aligned}
\]
The centers are real and distinct, so $\alpha_{ij}\ne0$ and
$\gamma_{ij}\ne0$. For $\{i,j,k\}=\{1,2,3\}$, it follows that
$t^2G_{ij}(t)G_{ik}(t)$ is an ordinary polynomial of degree four
with non-zero constant term.
Its four simple roots correspond to the intersection points identified
above. The polynomial $t^3S_i(t)$ has degree at most six and also
vanishes at all four roots. The factor theorem therefore gives
\[
 t^3S_i(t)=t^2G_{ij}(t)G_{ik}(t)R_i(t)
\]
for some $R_i\in\mathbb C[t]$ of degree at most two.
Write $R_i(t)=\rho_{i0}+\rho_{i1}t+\rho_{i2}t^2$. Then
$S_i(t)=G_{ij}(t)G_{ik}(t)\bigl(t^{-1}R_i(t)\bigr)$, with
\[
 t^{-1}R_i(t)=\rho_{i0}t^{-1}+\rho_{i1}+\rho_{i2}t.
\]
The functions $1,x_i(t),y_i(t)$ form a basis of
$\operatorname{span}_{\mathbb C}\{1,t,t^{-1}\}$, as follows from
\eqref{eq:complex-circle-parametrization} and $r_i>0$.
Hence there is a unique affine polynomial $A_i(x,y)$ such that
$A_i(x_i(t),y_i(t))=t^{-1}R_i(t)$.
For $t=r_ie^{\mathrm i\theta}$, both $S_i(t)$ and
$G_{ij}(t)G_{ik}(t)$ are real, and the latter is non-zero because
the real circles are pairwise disjoint. Thus $A_i$ is real-valued
on the real circle. Its imaginary part on $\mathbb R^2$ is an
affine polynomial vanishing on that circle; Lemma~\ref{lem:circle-restriction}
forces this polynomial to be zero. Hence $A_i$ has real coefficients.

On each circle $\vf$ now agrees with $\sum_i A_i\mathcal T_i$.
Lemma~\ref{lem:circle-restriction} therefore makes
$f_1f_2f_3$ divide both components of the difference
$\vf-\sum_{i=1}^{3}A_i\mathcal T_i$.
Since this difference has degree at most six and $f_1f_2f_3$ has degree six
we obtain
\[
 \vf=\sum_{i=1}^3 A_i\mathcal T_i+f_1f_2f_3h
\]
for a constant vector $h$. Comparing the degree-six terms and
cancelling $(x^2+y^2)^2$ gives
\[
 J(z)\,\sum_{i=1}^3\ell_i+(x^2+y^2)h=0.
\]
Taking its scalar product with $z$ gives $h=0$, then
$\sum_{i=1}^3\ell_i=0$.
The degree-five equation, after applying $-J$ and dividing by
$x^2+y^2$, is
\begin{equation}\label{eq:uniform-degree-five}
 (x^2+y^2)\Bigl[\Bigl(\sum_{i=1}^3 c_i\Bigr)z
       -\sum_{i=1}^3\ell_i C_i\Bigr]
       +2\Bigl[\sum_{i=1}^3\ell_i(C_i^Tz)\Bigr]z=0.
\end{equation}
Taking its scalar product with $z$, using $z^Tz=x^2+y^2$, and cancelling
the non-zero polynomial factor $x^2+y^2$ gives
$\sum_{i=1}^3\ell_i(C_i^Tz)
=-\bigl(\sum_{i=1}^3 c_i\bigr)(x^2+y^2)$.
Substituting this into \eqref{eq:uniform-degree-five} and cancelling
$x^2+y^2$ yields
\[
 \sum_{i=1}^3\ell_i C_i+\Bigl(\sum_{i=1}^3 c_i\Bigr)z=0.
\]
Taking the Jacobian matrix of this identity with respect to $z$
gives the second condition in \eqref{eq:uniform-affine-constraints}.
These same computations prove sufficiency of
\eqref{eq:uniform-affine-constraints} for the degree cancellations;
invariance follows from the factors and tangency of each $\mathcal T_i$.
Finally a zero representation restricts to $A_i\mathcal T_i=0$ on
each circle, forcing the affine polynomial $A_i$ to vanish identically.
This proves uniqueness.
\end{proof}

The no-common-point hypothesis holds in both cases. For collinear centers
with $\Delta\ne0$ \eqref{eq:pencil-relation} excludes a common finite
complex point. For non-collinear centers any such point must be the unique
real intersection of two independent radical axes, contradicting pairwise
disjointness of the real circles.
\par

\subsection{The non-coaxial case}

\begin{prop}\label{prop:collinear-normal-form}
Assume $\Delta\ne0$.  A real polynomial differential system of degree at
most four has all three circles \eqref{eq:collinear-circles} as invariant
algebraic curves if and only if
\begin{equation}\label{sys:collinear-normal-form}
 \quad \binom{P}{Q}=\lambda\mathcal Z_1+\ell\mathcal Z_0,\quad
\end{equation}
where $\lambda\in\mathbb R$ and
$\ell=\alpha x+\beta y+\gamma\in\mathbb R[x,y]$.
\end{prop}

\begin{proof}
By Lemma~\ref{lem:uniform-affine-representation} write
$\vf=\sum_{i=1}^3 A_i\mathcal T_i$, with $A_i=\ell_i+c_i$
satisfying \eqref{eq:uniform-affine-constraints}.
Since $C_1=(0,0)^T,C_2=(1,0)^T,C_3=(u,0)^T$ the matrix
$\sum_{i=1}^3 C_i(\nabla\ell_i)^T$ has zero second row,
so the constraints reduce to
\[
 \ell_1+\ell_2+\ell_3=0,\qquad \ell_2+u\ell_3=0,
 \qquad c_1+c_2+c_3=0.
\]
Set
\begin{equation}\label{eq:eta-nu}
 \eta=(\eta_1,\eta_2,\eta_3)=(1-u,u,-1),\qquad
 \nu=(\nu_1,\nu_2,\nu_3)=(1,0,-1).
\end{equation}
The first two constraints give $\ell_i=\eta_i(\alpha x+\beta y)$.
Since $u\ne0$ the vectors $\eta,\nu$ form a basis of the plane
of vectors whose components sum to zero. Hence
$(c_1,c_2,c_3)=\lambda\nu+\gamma\eta$, and
\[
 A_i=\lambda\nu_i+\eta_i(\alpha x+\beta y+\gamma).
\]
Substitution into the representation of $\vf$, using
\eqref{eq:Z-collinear}, gives \eqref{sys:collinear-normal-form}.
Conversely these coefficients satisfy
\eqref{eq:uniform-affine-constraints}, so
Lemma~\ref{lem:uniform-affine-representation} proves sufficiency.
\end{proof}

Direct calculation gives $\vf(f_i)=K_if_i$, $i=1,2,3$, with
\begin{equation}\label{eq:cofactors-collinear}
 \begin{split}
 K_1={}&2uy\{\lambda f_2+\ell(f_2-f_3)\},\\
 K_2={}&2y\{\lambda[f_3-(1-u)f_1]
              +(1-u)\ell(f_3-f_1)\},\\
 K_3={}&2uy\{\lambda f_2+(1-u)\ell(f_2-f_1)\}.
 \end{split}
\end{equation}
\par

\begin{prop}\label{prop:collinear-degree-gcd}
For the non-coaxial family \eqref{sys:collinear-normal-form}:
\begin{enumerate}
\item the system has exact degree four if and only if
\[
 (\lambda,\alpha,\beta)\ne(0,0,0);
\]
\item if $\lambda\ne0$ then $P,Q$ are relatively prime if and only if
$Q(x,0)\not\equiv0$ and
\begin{equation}\label{eq:collinear-gcd-degeneracies}
 (\alpha,\beta)\ne(0,0)
 \quad\text{or}\quad
 \gamma\notin\left\{0,-\lambda,-\frac{\lambda}{1-u}\right\}.
\end{equation}
When $\alpha=\beta=0$ the excluded values
$\gamma=-\lambda/(1-u),0,-\lambda$ correspond respectively to the common
factors $f_1,f_2,f_3$.
\item if $\lambda=0$ and the system has exact degree four then the
non-constant polynomial $\ell$ is a common factor.  Thus no such system is
coprime.
\end{enumerate}
\end{prop}

\begin{proof}
For the columns $\mathcal Z_1,x\mathcal Z_0,y\mathcal Z_0$ take the
coefficients of $y^4,x^3y$ in the first component and $x^4$ in the
second. The resulting minor is
\[
 \begin{vmatrix}
  0&0&-\Delta\\
  2u&-\Delta-4u(1-u)&0\\
  -u&\Delta+2u(1-u)&0
 \end{vmatrix}
 =-u\Delta^2\ne0.
\]
Thus their quartic homogeneous parts are linearly independent, proving
the first assertion.

Suppose $\lambda\ne0$, and let $h$ be a common irreducible factor of
$P,Q$ not proportional to any of $y,f_1,f_2,f_3$.
The invariance identities imply $h\mid K_i$ for $i=1,2,3$.
By \eqref{eq:cofactors-collinear} and \eqref{eq:pencil-relation}
\[
 K_1-K_3=-2uy\ell\Delta.
\]
Since $u\Delta\ne0$ and $h\nmid y$ we have $h\mid\ell$.
The formula for $K_1$ then gives $h\mid f_2$. Since $f_2$ is
irreducible this contradicts the choice of $h$.
Thus every common irreducible factor is proportional to one of
$y,f_1,f_2,f_3$.

Since $y\mid P$ it is a common factor exactly when
$Q(x,0)\equiv0$. In the representation
$\vf=\sum_{i=1}^3 A_i\mathcal T_i$ the coefficients are
\[
 A_1=\lambda+(1-u)\ell,\qquad A_2=u\ell,\qquad
 A_3=-(\lambda+\ell).
\]
On $f_i=0$ the vector field equals $A_i\mathcal T_i$, and
$\mathcal T_i$ is non-zero along the real circle.
Hence $f_i$ is a common factor exactly when $A_i$ vanishes on that
circle. Since $A_i$ is affine this is equivalent to $A_i\equiv0$.
These three possibilities occur precisely when $\alpha=\beta=0$ and
$\gamma=-\lambda/(1-u),0,-\lambda$, respectively. Excluding them gives
\eqref{eq:collinear-gcd-degeneracies}, proving the second assertion.

When $\lambda=0$ exact degree four requires $\ell$ to be non-constant,
and $\vf=\ell\mathcal Z_0$ proves the last assertion.
\end{proof}

\subsection{Periodic orbits, isolation and hyperbolicity}

Write $(s_1,s_2,s_3)=(0,1,u)$ for the center abscissae.
On circle $f_k=0$ use the counterclockwise angular coordinate
$\theta_k$ defined by
\[
 x=s_k+r_k\cos\theta_k,\qquad y=r_k\sin\theta_k,
 \qquad k=1,2,3.
\]
For the vector field \eqref{sys:collinear-normal-form} the angular velocities are
\begin{equation}\label{eq:angular-collinear}
 \dot\theta_1=[\lambda+(1-u)\ell]f_2f_3,\qquad
 \dot\theta_2=u\ell f_1f_3,\qquad
 \dot\theta_3=-(\lambda+\ell)f_1f_2.
\end{equation}
Here each right-hand side is evaluated on the corresponding circle.
By pairwise disjointness the factors $f_i$ appearing there do not vanish,
so each circle is a periodic orbit exactly when its angular velocity is
nowhere zero. On $f_k=0$ we have
$\ell=\alpha s_k+\gamma+r_k(\alpha\cos\theta_k+\beta\sin\theta_k)$,
whose oscillatory term has amplitude $r_k\sqrt{\alpha^2+\beta^2}$.
Thus all three circles are periodic orbits if and only if
\begin{equation}\label{eq:periodic-collinear}
 |\lambda+(1-u)\gamma|>|1-u|r_1
               \sqrt{\alpha^2+\beta^2},\ 
 |\alpha+\gamma|>r_2\sqrt{\alpha^2+\beta^2},\ 
 |\lambda+\alpha u+\gamma|>r_3
               \sqrt{\alpha^2+\beta^2}.
\end{equation}
\par

\begin{prop}\label{prop:collinear-limit-cycles}
Assume $\Delta\ne0$ and \eqref{eq:periodic-collinear} holds. All three circles
\eqref{eq:collinear-circles} are algebraic limit cycles of
\eqref{sys:collinear-normal-form} if and only if
\begin{equation}\label{eq:limit-cycles-collinear}
 \lambda\beta\ne0.
\end{equation}
When these conditions hold the system is a coprime exact quartic polynomial differential system and
all three limit cycles are hyperbolic.
\end{prop}

\begin{proof}
By \eqref{eq:periodic-collinear} it remains to determine when the periodic
circles are isolated. Using $\eta_k,\nu_k$ from \eqref{eq:eta-nu} put
\[
 a_k^\pm(\theta)=\lambda\nu_k+
 \eta_k\bigl[\alpha(s_k+r_k\cos\theta)+\gamma
             \pm\beta r_k\sin\theta\bigr],
\]
and denote by $F_{ik}(\theta)$ the restriction of $f_i$ to the circle
$f_k=0$. Define the cofactor integrand by
$\mathcal I_k(\theta)=\left.K_k\right|_{f_k=0}(\theta)/\dot\theta_k(\theta)$,
using \eqref{eq:cofactors-collinear} and \eqref{eq:angular-collinear}.
For $\{i,j,k\}=\{1,2,3\}$ pairing the values at
$\theta$ and $2\pi-\theta$ gives
\begin{equation}\label{eq:paired-collinear}
 \mathcal I_k(\theta)+\mathcal I_k(2\pi-\theta)
 =-\frac{4\lambda\beta u\Delta r_k^2\sin^2\theta}
 {F_{ik}(\theta)F_{jk}(\theta)a_k^+(\theta)a_k^-(\theta)}.
\end{equation}
Under disjointness and \eqref{eq:periodic-collinear} every factor in the
denominator is non-zero and has constant sign. If $\lambda\beta\ne0$
the right-hand side has a fixed non-zero sign for $0<\theta<\pi$, so
$\int_0^{2\pi}\mathcal I_k(\theta)\,d\theta\ne0$.
By \eqref{eq:preliminary-angle-exponent} the positive-time
characteristic exponent is also non-zero.
Thus all three periodic circles are hyperbolic limit cycles, proving sufficiency.

If $\beta=0$ the system is reversible with respect to
$y=0$, and each prescribed periodic circle meets this line transversely
at two points. Lemma~\ref{lem:reversible-period-annulus} excludes isolation.
If $\lambda=0$ then $\vf=\ell\mathcal Z_0$, with $\ell\ne0$
near each periodic circle. Since $\mathcal Z_0$ is reversible with
respect to $y=0$ the same lemma and invariance of trajectories under
non-vanishing scalar multiplication again exclude isolation.
This proves necessity.

Finally $\beta\ne0$ excludes the common factors $f_i$ by
\eqref{eq:collinear-gcd-degeneracies}, and periodicity excludes $y$,
because each circle intersects $y=0$.
Proposition~\ref{prop:collinear-degree-gcd} therefore gives coprimality and,
as $\lambda\ne0$, exact degree four.
\end{proof}

\subsection{The coaxial case}

\begin{prop}\label{prop:collinear-coaxial}
Assume $\Delta=0$.  Put
\begin{equation}\label{eq:Y-coaxial}
 \mathcal Y=\binom{y(f_1-f_2)}{(1-x)f_1+xf_2}.
\end{equation}
Then a polynomial differential system of degree at most four leaves the
three circles invariant if and only if
\begin{equation}\label{sys:collinear-coaxial}
 \binom{P}{Q}=q_2(x,y)\mathcal Y,
\end{equation}
where $q_2\in\mathbb R[x,y]$ is an arbitrary polynomial of $\deg q_2\le2$.
Such a system has exact degree four if
and only if $\deg q_2=2$, while coprimality forces $q_2$ to be a non-zero
constant and hence reduces the degree to two.  Moreover
\begin{equation}\label{eq:first-integral-coaxial}
 I=\frac{f_1}{f_2}
\end{equation}
is a rational first integral wherever it is defined.
None of the three prescribed circles is a limit cycle.
\end{prop}

\begin{proof}
Direct calculation gives $\mathcal Y(f_1)=2yf_1$ and
$\mathcal Y(f_2)=2yf_2$. The identity
\eqref{eq:pencil-relation} gives the same identity
for $f_3$, proving sufficiency and the first integral follows.

For necessity let $\vf=(P,Q)$ be a polynomial vector field of degree
at most four leaving the three prescribed circles invariant, and put
\[
 H=f_2\vf(f_1)-f_1\vf(f_2)
   =2(\mathcal Y_2P-\mathcal Y_1Q).
\]
Clearly $f_1f_2\mid H$.  The identity
$uH=f_3\vf(f_1)-f_1\vf(f_3)$ on $\Delta=0$ also gives $f_3\mid H$.
The three factors are relatively prime and $\deg H\le6$, so
$H=c f_1f_2f_3$ for a real constant $c$.

Set
\[
 a=\frac{1+r_1^2-r_2^2}{2},\qquad h^2=a^2-r_1^2,
 \qquad \xi=x-a.
\]
Strict disjointness gives
\[
 4h^2=\bigl(1-(r_1+r_2)^2\bigr)
       \bigl(1-(r_1-r_2)^2\bigr)>0.
\]
Take $h>0$. With the center abscissae denoted by
$(s_1,s_2,s_3)=(0,1,u)$, the identity
\eqref{eq:pencil-relation} gives
$r_i^2=(s_i-a)^2-h^2$, $i=1,2,3$, and hence
\[
 f_i=\xi^2+y^2+h^2-2(s_i-a)\xi,\qquad
 \mathcal Y=(2\xi y,-\xi^2+y^2+h^2).
\]
At $(\xi,y)=(h,0)$ the field $\mathcal Y$ vanishes, whereas
$f_i=2h(h-(s_i-a))\ne0$ for every $i$, because $r_i>0$.
Evaluating $H=c f_1f_2f_3$ there forces $c=0$.
Thus $\mathcal Y_2P=\mathcal Y_1Q$.
The polynomials $2\xi y$ and $-\xi^2+y^2+h^2$ are relatively
prime, so $\vf=q_2\mathcal Y$ for a unique polynomial $q_2$ with
$\deg q_2\le2$. Since $\deg\mathcal Y=2$ the assertions about
exact degree four and coprimality follow.

Finally the only real zeros of $\mathcal Y$ are
$(\xi,y)=(\pm h,0)$, neither of which lies on a prescribed circle.
Near circles 1 and 3 the first integral $f_1/f_2$ has non-zero gradient,
since $\nabla(f_1/f_2)=2(\mathcal Y_2,-\mathcal Y_1)^T/f_2^2$.
Near circle 2 use $f_2/f_1$, whose gradient is likewise non-zero.
If $q_2$ has no zero on a prescribed circle that circle
is periodic and Lemma~\ref{lem:regular-first-integral} places it in a
period annulus. Otherwise it contains an equilibrium. Hence none is a
limit cycle.
\end{proof}

The relation $r_i^2=(s_i-a)^2-h^2$, with $h>0$, places every center
abscissa in either $( -\infty,a-h)$ or $(a+h,\infty)$.
Circles whose centers lie in the same interval are nested, whereas those
whose centers lie in different intervals are exterior to one another.
Thus the coaxial case yields configuration 8 when all three centers lie in
the same interval, and configuration 6 otherwise.
\par

\begin{re}[Scope of the coaxial distinction]
\label{rem:coaxial-split-scope}
The distinction between coaxial and non-coaxial triples is needed only
for three distinct collinear centers. A concentric pair has no finite
radical axis, while a common radical axis for circles with distinct centers
forces their centers to lie on a line perpendicular to that axis.

Excluding the coaxial case does not exclude configurations 6 or 8:
the non-coaxial examples in \eqref{eq:collinear-examples} realize all four
configurations 5--8 with three circular limit cycles.
\end{re}

\subsection{Hyperbolic realizations}

To realize each configuration in Table~\ref{tab:collinear-domains} take
$\lambda=\beta=1$, $\alpha=0$, and choose the remaining parameters
as follows.
\begin{equation}\label{eq:collinear-examples}
\begin{array}{c|c|c|c|c|c}
\text{configuration}&u&r_1&r_2&r_3&\gamma\\ \hline
5&3&1/4&1/5&1/6&1\\
6&4&2&1/2&1/3&4\\
7&-1&4&1/2&1/2&4\\
8&1/2&4&2&1/4&4
\end{array}
\end{equation}
Each row satisfies the corresponding geometric inequalities, $\Delta\ne0$,
and \eqref{eq:periodic-collinear};
Proposition~\ref{prop:collinear-limit-cycles} therefore gives three hyperbolic
circular limit cycles.

Together with Propositions~\ref{prop:collinear-normal-form},
\ref{prop:collinear-limit-cycles} and~\ref{prop:collinear-coaxial} these
examples complete the proof of Theorem~\ref{main:collinear}.
\par

\section{Three non-collinear centers: configurations 9, 10, 11 and 12}
\label{sec:noncollinear}

This section proves Theorem~\ref{main:noncollinear}.

A translation, a rotation, a homothety and, if necessary, a reflection put the centers at
\begin{equation}\label{eq:noncollinear-circles}
 C_1=(0,0),\qquad C_2=(1,0),\qquad C_3=(u,v),
\end{equation}
with $u\in\mathbb{R}$ and $v>0$. The corresponding circle polynomials are
\begin{equation}\label{eq:noncollinear-fi}
 f_1=x^2+y^2-r_1^2,\qquad
 f_2=(x-1)^2+y^2-r_2^2,\qquad
 f_3=(x-u)^2+(y-v)^2-r_3^2,
\end{equation}
where $r_i>0$ for $i=1,2,3$. Put
\begin{equation}\label{eq:noncollinear-distances}
 d_{12}=1,\qquad d_{13}=\sqrt{u^2+v^2},\qquad
 d_{23}=\sqrt{(u-1)^2+v^2}.
\end{equation}
Throughout this section the real circles are also assumed to be pairwise disjoint: for each $i<j$ either $d_{ij}>r_i+r_j$ or
$d_{ij}<|r_i-r_j|$. Their four possible containment types are specified
in Table~\ref{tab:noncollinear-domains}.
The disk labels follow the containment roles indicated in each row.

\begin{table}[htbp]
\centering
\renewcommand{\arraystretch}{1.3}
\begin{tabular}{c|l|l}
Configuration & Geometry & Strict inequalities \\ \hline
9 & Pairwise exterior disks
  & $r_1+r_2<1,\ r_1+r_3<d_{13},$ \\
  & & $r_2+r_3<d_{23}$ \\
10 & $D_2\subset\operatorname{int}D_1$, $D_3$ exterior to $D_1$
  & $1+r_2<r_1,\ r_1+r_3<d_{13}$ \\
11 & $D_2,D_3\subset\operatorname{int}D_1$, mutually exterior
  & $1+r_2<r_1,\ d_{13}+r_3<r_1,$ \\
  & & $r_2+r_3<d_{23}$ \\
12 & $D_3\subset\operatorname{int}D_2\subset\operatorname{int}D_1$
  & $1+r_2<r_1,\ d_{23}+r_3<r_2$ \\
\end{tabular}
\par\medskip
\caption{Domains for configurations 9--12, with $v>0$ and
$r_1,r_2,r_3>0$.}
\label{tab:noncollinear-domains}
\end{table}

The other exterior or containment inequalities in configurations 10 and 12
follow from the inequalities listed in Table~\ref{tab:noncollinear-domains}
and the triangle inequality.

We now derive a normal form valid for all four configurations.
Introduce the degree-five fields tangent to the individual circles,
\begin{equation}\label{eq:T-noncollinear}
 \mathcal T_1=f_2f_3\binom{-y}{x},\qquad
 \mathcal T_2=f_1f_3\binom{-y}{x-1},\qquad
 \mathcal T_3=f_1f_2\binom{-(y-v)}{x-u}.
\end{equation}
Define the affine polynomials
\begin{equation}\label{eq:L-noncollinear}
 L_1=-vx+(u-1)y-v,\qquad L_2=vx-uy,\qquad L_3=y.
\end{equation}
For the non-collinear case define the three quartic fields below, using
$\mathcal Z_0,\mathcal Z_1,\mathcal Z_2$ as notation specific to this case.
\begin{equation}\label{eq:Z-noncollinear}
 \mathcal Z_1=\mathcal T_2-\mathcal T_1,\qquad
 \mathcal Z_2=\mathcal T_3-\mathcal T_1,\qquad
 \mathcal Z_0=L_1\mathcal T_1+L_2\mathcal T_2+L_3\mathcal T_3.
\end{equation}
The degree-six and degree-five parts cancel in $\mathcal Z_0$, while the
degree-five parts cancel in $\mathcal Z_1,\mathcal Z_2$.  Direct comparison
of their quartic homogeneous parts shows that these three fields are linearly
independent whenever $v\ne0$.

\begin{prop}\label{prop:noncollinear-normal-form}
Assume that the circles \eqref{eq:noncollinear-fi} are pairwise
disjoint and that $v>0$.
A real polynomial differential system of degree at most four has the three
circles \eqref{eq:noncollinear-fi}, whose centers are non-collinear,
as invariant algebraic curves
if and only if
\begin{equation}\label{sys:noncollinear-normal-form}
 \quad \binom{P}{Q}=\lambda\mathcal Z_1+
 \mu\mathcal Z_2+\kappa\mathcal Z_0,\quad
\end{equation}
where $\lambda,\mu,\kappa\in\R$.  Equivalently
\begin{equation}\label{eq:A-noncollinear}
 \binom{P}{Q}=A_1\mathcal T_1+A_2\mathcal T_2+A_3\mathcal T_3,
\end{equation}
where
\begin{equation}\label{eq:A123-noncollinear}
 A_1=-\lambda-\mu+\kappa L_1,\qquad
 A_2=\lambda+\kappa L_2,\qquad
 A_3=\mu+\kappa L_3.
\end{equation}
If $(\lambda,\mu,\kappa)\ne(0,0,0)$ the vector field in
\eqref{sys:noncollinear-normal-form} has degree exactly four.
\end{prop}

\begin{proof}
The two independent real radical-axis equations force any common finite
complex point of the three circles to be real, which strict disjointness
excludes. Thus Lemma~\ref{lem:uniform-affine-representation} applies. Write
its affine coefficients as $A_i=\ell_i+c_i$, where
$\ell_i=a_{i1}x+a_{i2}y$ and $c_i\in\mathbb R$.
For the normalized centers
\[
 \sum_i C_i(\nabla\ell_i)^T=\begin{pmatrix}
 a_{21}+ua_{31}&a_{22}+ua_{32}\\
 va_{31}&va_{32}
 \end{pmatrix}=-\Bigl(\sum_i c_i\Bigr)\mathrm{Id}_2.
\]
Set $\kappa=-(\sum_i c_i)/v$. Since $v>0$ these equations
together with $\sum_i\ell_i=0$ uniquely determine the following expressions:
\[
 \ell_1=\kappa(-vx+(u-1)y),\quad
 \ell_2=\kappa(vx-uy),\quad \ell_3=\kappa y.
\]
Since $\sum_i c_i=-\kappa v$ we may put
$c_2=\lambda,c_3=\mu,c_1=-\lambda-\mu-\kappa v$.
This is precisely \eqref{eq:A123-noncollinear}. Conversely these
coefficients satisfy the constraints in
Lemma~\ref{lem:uniform-affine-representation} and therefore define a vector
field of degree at most four leaving all three circles invariant.

For exact degree form the coefficient matrix of the quartic homogeneous
parts of $P,Q$, with columns ordered as $(\lambda,\mu,\kappa)$.
The rows for $y^4,xy^3,x^3y$ in $P$, and those for the same monomials
in $Q$, give respectively the two $3\times3$ minors
\[
 -8(2u-1)v^2,\qquad -8v(u^2-u-v^2).
\]
For real $v\ne0$ they cannot both vanish: the first would force
$u=1/2$, where the bracket in the second is $-1/4-v^2\ne0$.
Thus a vanishing quartic part forces $\lambda=\mu=\kappa=0$.
\end{proof}

The corresponding cofactors are
\begin{equation}\label{eq:cofactors-noncollinear}
 \begin{aligned}
 K_1&=2\bigl[-A_2f_3y+A_3f_2(vx-uy)\bigr],\\
 K_2&=2\bigl[A_1f_3y+A_3f_1\bigl(v(x-1)-(u-1)y\bigr)\bigr],\\
 K_3&=2\bigl[A_1f_2(uy-vx)+A_2f_1\bigl((u-1)y-v(x-1)\bigr)\bigr].
 \end{aligned}
\end{equation}
\par
\subsection{Periodic orbits, coprimality and hyperbolicity}
\label{subsec:noncollinear-periodicity}

On the $k$-th circle use the counterclockwise parametrization
\begin{equation}\label{eq:Gamma-noncollinear}
 \Gamma_k(\theta)=C_k+r_k(\cos\theta,\sin\theta),\qquad k=1,2,3.
\end{equation}
Since the three circles are pairwise disjoint the restrictions of the two
other circle factors have constant non-zero signs.  Equations
\eqref{eq:A-noncollinear} and \eqref{eq:Gamma-noncollinear} give
\begin{equation}\label{eq:angular-noncollinear}
 \left.\dot\theta_k\right|_{f_k=0}=A_k f_i f_j,
 \qquad \{i,j,k\}=\{1,2,3\}.
\end{equation}
Here $A_k$ is defined in \eqref{eq:A123-noncollinear}, and all
polynomials on the right-hand side are evaluated on the $k$-th circle
at $\Gamma_k(\theta)$.
On the $k$-th circle $A_k$ has the form
$a+b\cos\theta+c\sin\theta$.
Therefore all three invariant circles are periodic orbits if and only if
\begin{equation}\label{eq:periodic-noncollinear}
 |\lambda+\mu+\kappa v|>|\kappa|r_1d_{23},\qquad
 |\lambda+\kappa v|>|\kappa|r_2d_{13},\qquad
 |\mu+\kappa v|>|\kappa|r_3.
\end{equation}
\par

Under \eqref{eq:periodic-noncollinear}, let $\chi_k$
denote the characteristic exponent of the periodic circle $f_k=0$,
for $k=1,2,3$, and define
\begin{equation}\label{eq:chi-noncollinear}
 I_k=\int_0^{2\pi}
 \frac{K_k(\Gamma_k(\theta))}
 {A_k(\Gamma_k(\theta))
  f_i(\Gamma_k(\theta))f_j(\Gamma_k(\theta))}\,d\theta,
 \qquad \{i,j,k\}=\{1,2,3\}.
\end{equation}
Equation~\eqref{eq:preliminary-angle-exponent} then gives
\begin{equation}\label{eq:time-exponent-noncollinear}
 \chi_k=\varepsilon_k I_k,\qquad
 \varepsilon_k=\operatorname{sgn}\!\left(\left.\dot\theta_k\right|_{f_k=0}\right).
\end{equation}
Here $\varepsilon_k\in\{-1,1\}$ is constant, so $\chi_k=0$ if and
only if $I_k=0$.
\par

We now evaluate $I_1,I_2,I_3$ explicitly.

\begin{prop}\label{prop:noncollinear-hyperbolic}
Assume that the circles \eqref{eq:noncollinear-fi} are pairwise disjoint
and that \eqref{eq:periodic-noncollinear} holds. Then the vector field in
\eqref{sys:noncollinear-normal-form} is coprime and has degree exactly four.
The quantities $\delta_k,H_{ki},B_{ki}$, and $\Phi_k$
are as defined before Theorem~\ref{main:noncollinear}, and
$B_{ik}=B_{ki}$ for $i\ne k$.
Moreover $\delta_k\ne0$ and
$H_{ki}\ne0$ for $k\ne i$, and
\begin{equation}\label{eq:explicit-exponents-noncollinear}
 I_k=\frac{2\pi\kappa r_k^2}{\delta_k}
       \left(\frac{B_{ki}}{H_{ki}}+\frac{B_{kj}}{H_{kj}}\right)
     =\frac{2\pi\kappa r_k^2\Phi_k}
            {\delta_kH_{ki}H_{kj}},
 \qquad \{i,j,k\}=\{1,2,3\}.
\end{equation}
Consequently the three circles are simultaneously hyperbolic algebraic
limit cycles if and only if
\begin{equation}\label{eq:chi-product-noncollinear}
 \kappa\Phi_1\Phi_2\Phi_3\ne0.
\end{equation}
\end{prop}

\begin{proof}
We first prove exact degree four and coprimality using only the periodicity
conditions \eqref{eq:periodic-noncollinear}, without assuming hyperbolicity.
Periodicity excludes the zero vector field, so exact degree four follows
from Proposition~\ref{prop:noncollinear-normal-form}.
If a non-constant polynomial $h$ divided both components periodicity
would exclude $f_i\mid h$. Since each positive-radius circle polynomial
is irreducible we have $\gcd(f_i,h)=1$. Write $\vf=h\mathcal Y$. The
identity $h\mathcal Y(f_i)=K_if_i$ then implies
$f_i\mid\mathcal Y(f_i)$. Thus the non-zero field $\mathcal Y$
would still leave all three circles invariant and have degree at most
three, contradicting Proposition~\ref{prop:noncollinear-normal-form}.

For a uniform notation write $C_i=(a_i,b_i)$, with
$(a_1,b_1)=(0,0)$, $(a_2,b_2)=(1,0)$, and
$(a_3,b_3)=(u,v)$, and set
\begin{equation}\label{eq:Dki-noncollinear}
 D_{ki}(x,y)=-(x-a_k)(y-b_i)+(y-b_k)(x-a_i).
\end{equation}
The cofactors in \eqref{eq:cofactors-noncollinear} can then be written as
\[
 K_k=2\bigl(A_if_jD_{ki}+A_jf_iD_{kj}\bigr),
 \qquad \{i,j,k\}=\{1,2,3\}.
\]
Differentiation along $\Gamma_k$ gives
$(f_i\circ\Gamma_k)'=-2D_{ki}\circ\Gamma_k$.
Combining these identities rewrites
\eqref{eq:chi-noncollinear} as
\[
 I_k=-\sum_{i\ne k}\int_0^{2\pi}
 \frac{A_i(\Gamma_k(\theta))}{A_k(\Gamma_k(\theta))}
 \frac{(f_i\circ\Gamma_k)'(\theta)}{f_i(\Gamma_k(\theta))}\,d\theta.
\]

To evaluate these integrals, we first derive a general integration formula.
Write $\boldsymbol n(\theta)=(\cos\theta,\sin\theta)$ and let
\[
 L=a+\boldsymbol b\cdot\boldsymbol n,\qquad M=c+\boldsymbol d\cdot\boldsymbol n,\qquad N=e+\boldsymbol g\cdot\boldsymbol n,
\]
where $a,c,e\in\mathbb R$ are scalars and
\[
 \boldsymbol b=(b_x,b_y),\qquad
 \boldsymbol d=(d_x,d_y),\qquad
 \boldsymbol g=(g_x,g_y)
\]
are vectors in $\mathbb R^2$. Assume
$|a|>\lVert\boldsymbol b\rVert$ and
$|c|>\lVert\boldsymbol d\rVert$, where $\lVert\cdot\rVert$ denotes the
Euclidean norm. We write
$\det(\boldsymbol b,\boldsymbol d)=b_xd_y-b_yd_x$ for the determinant
with columns $\boldsymbol b$ and $\boldsymbol d$. Put
\[
 \delta=\operatorname{sgn}(a)\sqrt{a^2-\lVert\boldsymbol b\rVert^2},\qquad
 \eta=\operatorname{sgn}(c)\sqrt{c^2-\lVert\boldsymbol d\rVert^2},\qquad
 \mathcal H=ac+\delta\eta-\boldsymbol b\cdot \boldsymbol d.
\]
We claim that $\mathcal H\ne0$ and
\begin{equation}\label{eq:affine-trigonometric-integral}
 -\int_0^{2\pi}\frac{NM'}{LM}\,d\theta
 =\frac{2\pi\det(\boldsymbol b,\boldsymbol g)}{\delta(a+\delta)}
  +\frac{2\pi\bigl(e\det(\boldsymbol b,\boldsymbol d)-a\det(\boldsymbol g,\boldsymbol d)-c\det(\boldsymbol b,\boldsymbol g)\bigr)}
         {\delta\mathcal H},
\end{equation}
where the prime denotes differentiation with respect to $\theta$.

To prove the claim identify $\boldsymbol b,\boldsymbol d$ with the complex numbers
$b_x+\mathrm i b_y,d_x+\mathrm i d_y$, and set
\[
 p=-\frac{b_x+\mathrm i b_y}{a+\delta},\qquad
 q=-\frac{d_x+\mathrm i d_y}{c+\eta}.
\]
Then $|p|,|q|<1$, and direct expansion gives
\[
 L=\frac{a+\delta}{2}|e^{\mathrm i\theta}-p|^2,\qquad
 M=\frac{c+\eta}{2}|e^{\mathrm i\theta}-q|^2,
\]
and
\begin{equation}\label{eq:nonzero-exponent-denominator}
 2\mathcal H=(a+\delta)(c+\eta)|1-q\overline p|^2\ne0.
\end{equation}
We use Fourier orthogonality and geometric summation to evaluate the
integrals explicitly. Since $|p|,|q|<1$ the following Fourier series
converge absolutely and uniformly, justifying termwise integration:
\[
 \frac1L=\frac1\delta
 \left[1+\sum_{m=1}^{\infty}
 \bigl(p^m e^{-\mathrm i m\theta}
       +\overline p^{\,m}e^{\mathrm i m\theta}\bigr)\right],\qquad
 \frac{M'}M=\mathrm i\sum_{m=1}^{\infty}
 \bigl(q^m e^{-\mathrm i m\theta}
       -\overline q^{\,m}e^{\mathrm i m\theta}\bigr).
\]
The first gives
$\int_0^{2\pi}e^{\mathrm i m\theta}/L\,d\theta=2\pi p^m/\delta$
for every integer $m\ge0$. Using these moments and the identity
$M'=-d_x\sin\theta+d_y\cos\theta$ we evaluate the integral of $M'/L$.
For the integral of $M'/(LM)$ set $z=q\overline p$, so
$|z|<1$. Integrating the second Fourier series against $1/L$ gives
\[
 J_2:=\int_0^{2\pi}\frac{M'}{LM}\,d\theta
 =\frac{2\pi\mathrm i}{\delta}
   \sum_{m=1}^{\infty}\bigl(z^m-\overline z^{\,m}\bigr)
 =-\frac{4\pi}{\delta}\operatorname{Im}\frac{z}{1-z}.
\]
The identities
\[
 \operatorname{Im}z=
 \frac{\det(\boldsymbol b,\boldsymbol d)}{(a+\delta)(c+\eta)},\qquad
 \operatorname{Im}\frac{z}{1-z}
 =\frac{\operatorname{Im}z}{|1-z|^2}
\]
and \eqref{eq:nonzero-exponent-denominator} give the value of
$J_2$. Together the two calculations yield
\[
 \int_0^{2\pi}\frac{M'}L\,d\theta
 =-\frac{2\pi\det(\boldsymbol b,\boldsymbol d)}{\delta(a+\delta)},\qquad
 \int_0^{2\pi}\frac{M'}{LM}\,d\theta
 =-\frac{2\pi\det(\boldsymbol b,\boldsymbol d)}{\delta\mathcal H}.
\]

If $\det(\boldsymbol b,\boldsymbol d)\ne0$ write $N=\alpha L+\beta M+\gamma$, with
\[
 \alpha=\frac{\det(\boldsymbol g,\boldsymbol d)}{\det(\boldsymbol b,\boldsymbol d)},\qquad
 \beta=\frac{\det(\boldsymbol b,\boldsymbol g)}{\det(\boldsymbol b,\boldsymbol d)},\qquad
 \gamma=e-\alpha a-\beta c.
\]
Since $M$ is periodic and nowhere zero we have
\[
 \int_0^{2\pi}\frac{M'}M\,d\theta
 =\bigl[\log|M(\theta)|\bigr]_0^{2\pi}=0.
\]
Substituting the decomposition of $N$ and the two integral identities
therefore proves \eqref{eq:affine-trigonometric-integral}.
The strict inequalities $|a|>\lVert\boldsymbol b\rVert$ and
$|c|>\lVert\boldsymbol d\rVert$ keep $|L|$ and $|M|$ uniformly
bounded away from zero in a neighborhood of each parameter point, so the integral
depends continuously on the parameters. The right-hand side is also
continuous, since $\delta(a+\delta)\ne0$ and
\eqref{eq:nonzero-exponent-denominator} gives $\mathcal H\ne0$.
Parameters with $\det(\boldsymbol b,\boldsymbol d)\ne0$ are dense in
this domain; hence the identity extends to
$\det(\boldsymbol b,\boldsymbol d)=0$ by continuity.

For each ordered pair $(k,i)$ with $i\ne k$, let $j$
be the remaining index in $\{1,2,3\}$. We now apply
\eqref{eq:affine-trigonometric-integral} with
\[
 L=A_k\circ\Gamma_k,\qquad M=f_i\circ\Gamma_k,\qquad
 N=A_i\circ\Gamma_k.
\]
Their coefficients are
\[
\begin{aligned}
 a&=\omega_k,& \boldsymbol b&=r_k\nabla A_k,\\
 c&=c_{ki},& \boldsymbol d&=2r_k(C_k-C_i),\\
 e&=A_i(C_k),& \boldsymbol g&=r_k\nabla A_i.
\end{aligned}
\]
Indeed $A_k(C_k)=\omega_k$ and
$\lVert\boldsymbol b\rVert^2=\kappa^2r_k^2d_{ij}^2$, so periodicity gives $|a|>\lVert\boldsymbol b\rVert$
and $\delta=\delta_k\ne0$.
Disjointness gives
\[
 c_{ki}^2-\lVert\boldsymbol d\rVert^2
 =\bigl(d_{ki}^2-(r_k+r_i)^2\bigr)
  \bigl(d_{ki}^2-(r_k-r_i)^2\bigr)>0,
\]
hence $|c|>\lVert\boldsymbol d\rVert$ and $\eta=h_{ki}$.
The explicit gradients are
\[
 \nabla A_1=\kappa(-v,u-1),\qquad
 \nabla A_2=\kappa(v,-u),\qquad
 \nabla A_3=\kappa(0,1).
\]
For each of the six ordered pairs $(k,i)$ these give
\[
 \boldsymbol b\cdot \boldsymbol d=2\kappa vr_k^2,\qquad
 e\det(\boldsymbol b,\boldsymbol d)-a\det(\boldsymbol g,\boldsymbol d)-c\det(\boldsymbol b,\boldsymbol g)=\kappa r_k^2B_{ki}.
\]
Thus $\mathcal H=H_{ki}$, which is non-zero by
\eqref{eq:nonzero-exponent-denominator}.

For each fixed $k$ sum the two contributions to $I_k$.
Since $\nabla A_1+\nabla A_2+\nabla A_3=0$ the first terms on the
right-hand side of \eqref{eq:affine-trigonometric-integral} cancel:
\[
 \sum_{i\ne k}\det(r_k\nabla A_k,r_k\nabla A_i)
 =\det(r_k\nabla A_k,-r_k\nabla A_k)=0.
\]
The remaining terms give the first equality in
\eqref{eq:explicit-exponents-noncollinear}; combining the two fractions
gives the second equality by the definitions of $\Phi_k$.

Finally $\chi_k=\varepsilon_k I_k$, with
$\varepsilon_k\in\{-1,1\}$, and all denominators in
\eqref{eq:explicit-exponents-noncollinear} are non-zero; thus
\eqref{eq:chi-product-noncollinear} is equivalent to
$\chi_1\chi_2\chi_3\ne0$, proving the hyperbolicity criterion by
Section~\ref{sec:preliminary}.
\end{proof}
\par

\begin{prop}[The case $\kappa=0$]\label{prop:noncollinear-integrable}
For the vector field in \eqref{sys:noncollinear-normal-form} assume
$\kappa=0$. Then each of the circles \eqref{eq:noncollinear-fi}
either consists entirely of equilibria or belongs to a period annulus.
In particular none of the three circles is a limit cycle.
\end{prop}

\begin{proof}
When $\kappa=0$ the functions $A_i$ are constants:
$A_1=-\lambda-\mu$, $A_2=\lambda$, and $A_3=\mu$.
If $A_k=0$ then \eqref{eq:A-noncollinear} shows that every
point of $f_k=0$ is an equilibrium.

Now suppose $A_k\ne0$. Pairwise disjointness ensures
that the angular velocity $A_kf_if_j$, with
$\{i,j,k\}=\{1,2,3\}$, is nowhere zero on $f_k=0$.
Thus this circle is a periodic orbit.

On the complement of the circles
\begin{equation}\label{eq:first-integral-noncollinear}
 H=A_1\log|f_1|+A_2\log|f_2|+A_3\log|f_3|
\end{equation}
is a first integral because
\[
 \vf=\frac{f_1f_2f_3}{2}J(\nabla H),
\]
where $J(a,b)=(-b,a)$.
To remove the logarithmic singularity of $H$ on the $k$-th circle,
choose a neighborhood of this circle that does not meet the other two
circles and define there
\[
 G_k=f_k\prod_{i\ne k}|f_i|^{A_i/A_k}.
\]
The other circle polynomials do not vanish in this neighborhood, so
$G_k$ is smooth there. Away from $f_k=0$, the identity
$G_k=\operatorname{sgn}(f_k)\exp(H/A_k)$ gives $\vf(G_k)=0$;
by continuity this also holds on the circle. Moreover
\[
 \left.\nabla G_k\right|_{f_k=0}
 =\left(\prod_{i\ne k}|f_i|^{A_i/A_k}\right)\nabla f_k\ne0.
\]
Lemma~\ref{lem:regular-first-integral} therefore places this periodic
circle in a period annulus.
\end{proof}

\par

\subsection{Hyperbolic realizations}
\label{subsec:noncollinear-domains}

With $(\lambda,\mu,\kappa)=(10,10,1)$ the following examples
realize each of the four configurations in
Table~\ref{tab:noncollinear-domains} with three hyperbolic circular
limit cycles:
\begin{equation}\label{eq:noncollinear-examples}
\begin{array}{c|c|c|c|c|c|c|c|c}
\text{configuration}&u&v&r_1&r_2&r_3&I_1&I_2&I_3\\ \hline
9&3&2&1/5&1/4&3/20&-0.018845&-0.181265&-0.002775\\
10&4&2&2&1/2&1/2&-0.112088&0.076185&0.007041\\
11&1/2&1/2&2&1/5&1/5&-0.076336&0.038474&0.018228\\
12&1/2&1/5&3&1&1/5&-0.143075&0.604852&-0.028668
\end{array}
\end{equation}
Each row of \eqref{eq:noncollinear-examples} satisfies the corresponding
geometric inequalities in Table~\ref{tab:noncollinear-domains} and the
periodicity inequalities \eqref{eq:periodic-noncollinear} exactly.
The last three columns give decimal approximations to $I_1,I_2,I_3$.
To certify non-vanishing enclose every square root in
\eqref{eq:explicit-exponents-noncollinear} between consecutive multiples
of $10^{-12}$, checking the endpoints by squaring. Rational interval
arithmetic then gives $|I_k|/\pi>1/1500$ for all twelve entries, and
hence $|I_k|>1/500$.
Thus all three circles in each example are hyperbolic limit cycles by
Proposition~\ref{prop:noncollinear-hyperbolic}.

Propositions~\ref{prop:noncollinear-normal-form},
\ref{prop:noncollinear-hyperbolic}, and~\ref{prop:noncollinear-integrable}, together with the realizations in
\eqref{eq:noncollinear-examples}, complete the proof of
Theorem~\ref{main:noncollinear}.
\par

\raggedbottom
\section*{Data availability statement}

This paper has no data.

\section*{Author contribution statement}

Both authors contributed equally to this paper.

\section*{Funding declaration}

The first author is partially supported by the China Scholarship Council
(No. 202506230080) and the National Natural Science Foundation of China
(No. 12271355).  The second author is partially supported by the Agencia
Estatal de Investigaci\'{o}n of Spain grant PID2025-168432NB-I00.

\end{document}